\documentclass[11pt]{article}
\usepackage{amsmath,amstext,amsthm,amssymb,amsfonts,latexsym,amscd,bezier}
\usepackage[brazil,english]{babel}
\usepackage[latin1]{inputenc}
\usepackage{color}
\usepackage{hyperref}
\usepackage{pdfsync}
\newcommand{\R}{\Bbb{R}}

\newcommand{\s}{\Bbb{S}}

\newcommand{\co}{\Bbb{C}}
\newcommand{\C}{\mathbb{C}}
\newcommand{\He}{\mathbb{H}}

\newcommand{\normc}[3]{\| #1 \|^{#2}_{#3}}

\newtheorem{lema}{Lemma}[section]

\newtheorem{teorema}{Theorem}[section]
\newtheorem{coro}{Corollary}[section]
\newtheorem{propo}{Proposition}[section]
\newtheorem{conj}{Conjecture}[section]
\newcommand{\lp}{\left(}
\newcommand{\rp}{\right)}
\newcommand{\n}{\noindent}
\newcommand{\vs}{\vspace}

\newcommand{\al}{\alpha}
\newcommand{\oal}{\bar{\alpha}}
\newcommand{\be}{\beta}
\newcommand{\obe}{\bar{\beta}}
\newcommand{\ga}{\gamma}
\newcommand{\oga}{\bar{\gamma}}
\newcommand{\rro}{\rho}

\newcommand{\omi}{\overline{\mu}}
\newcommand{\om}{\omega}

\newcommand{\Ga}{\Gamma}

\newcommand{\ta}{\theta}
\newcommand{\tta}{\widetilde{\theta}}
\newcommand{\tih}{\widetilde{h}}
\newcommand{\Ta}{\Theta}
\newcommand{\la}{\lambda}
\newcommand{\La}{\Lambda}
\newcommand{\osi}{\bar{\sigma}}
\newcommand{\si}{\sigma}

\newcommand{\dl}{\delta}

\newcommand{\lap}{\Delta}

\newcommand{\grad}{\nabla}

\newcommand{\til}[1]{\widetilde{#1}}

\begin{document}

\title{The $k$-Yamabe problem on CR manifolds
\footnote{2010 Mathematics Subject Classification: 51Mxx; 53A30; 58Jxx}
\footnote{Key words: Conformal geometry, CR geometry, Yamabe problem, fully nonlinear equations}
}

\author{\textbf{Ezequiel Barbosa  \footnote{\textit{E-mail addresses}:
ezequiel@mat.ufmg.br (E. R. Barbosa)}}\\ {\small\it Departamento de Matem\'{a}tica, Universidade Federal de Minas Gerais,}\\ {\small\it Caixa Postal 702, 30123-970, Belo Horizonte, MG, Brazil}\\
\textbf{Luiz Gustavo Carneiro \footnote{\textit{E-mail addresses}:
luizgustavo@iceb.ufop.br (L. G. Carneiro)}}\\ {\small\it Departamento de Matem\'atica, Universidade Federal de Ouro Preto,}\\ {\small\it 35400-000, Ouro Preto, MG, Brazil}\\
\textbf{Marcos Montenegro \footnote{\textit{E-mail addresses}:
montene@mat.ufmg.br (M. Montenegro)}}\\ {\small\it Departamento de Matem\'{a}tica, Universidade Federal de Minas Gerais,}\\ {\small\it
Caixa Postal 702, 30123-970, Belo Horizonte, MG, Brazil}}

\date{}

\maketitle

\markboth{abstract}{abstract}
\addcontentsline{toc}{chapter}{abstract}

\hrule \vspace{0,2cm}

\n {\bf Abstract}

We introduce the notion of pseudohermitian $k$-curvature, which is a natural extension of the Webster scalar curvature, on an orientable manifold endowed with a strictly pseudoconvex pseudohermitian structure (referred here as a CR manifold) and raise the $k$-Yamabe problem on a compact CR manifold. When $k=1$, the problem was proposed and partially solved by Jerison and Lee for CR manifolds non-locally CR-equivalent to the CR sphere. For $k > 1$, the problem can be translated in terms of the study of a fully nonlinear equation of type complex $k$-Hessian. We provide some partial answers related to the CR $k$-Yamabe problem. We establish that its solutions with null Cotton tensor are critical points of a suitable geometric functional constrained to pseudohermitian structures of unit volume. Thanks to this variational property, we establish a Obata type result for the problem and also compute the infimum of the functional on the CR sphere. Furthermore, we show that this value is an upper bound for the corresponding one on any compact CR manifolds and, assuming the CR Yamabe invariant is positive, we prove that such an upper bound is only attained for compact CR manifolds locally CR-equivalent to the CR sphere. In the Riemannian field, recent advances have been produced in a series of outstanding works.

\vspace{0.5cm}
\hrule\vspace{0.2cm}

\n {\small {\bf Content}}

\n {\small 1. Introduction, overview and main results ............................................................................................................................\ 2\\
2. Pseudohermitian geometric tensors and the CR $k$-curvature ...........................................................................................\ 10\\
3. The CR $k$-Yamabe problem .............................................................................................................................................\ 16\\
4. The variational property and pseudo-Einstein manifolds .................................................................................................\ 18\\
5. An ellipticity property .....................................................................................................................................................\ 23\\
6. The CR $k$-Yamabe constant .............................................................................................................................................\ 24\\
References
............................................................................................................................................................................\ 30}

\vspace{0.5cm}
\hrule\vspace{0.2cm}

\pagebreak

\section{Introduction, overview and main statements}

The geometry of CR manifolds, namely the abstract model of real hypersurfaces in complex manifolds, has attracted, since the late 70's,
a lot of attention of prominent mathematicians as for instance Chern, Moser \cite{CM}, Fefferman \cite{Fe}, Jacobowitz \cite{Jac}, Jerison,
Lee \cite{JL1}, \cite{JL2}, \cite{JL3}, \cite{JL4}, \cite{Lee1}, \cite{Lee2}, Tanaka \cite{Tan}, Webster \cite{Web}, among many others. This
geometry is very rich when the CR manifold admits a strictly pseudoconvex structure in which case we encounter a great analogy with the geometry
of Riemannian manifolds.

A pseudohermitian structure on a CR manifold $(M, T_{1,0}(M))$ of dimension $2n + 1$, $n \geq 1$, is a
contact form $\theta$ that annihilates the Levi distribution $H(M)$
defined as the real part of the subspace $T_{1,0}(M) \oplus
T_{0,1}(M)$ of the complexified of $T(M)$. Such a structure determines a natural hermitian metric
on the CR structure $T_{1,0}(M)$, namely the Levi form $L_\theta$, which is well-defined modulo a multiplication by a smooth function,
exactly as occur in conformal Riemannian geometry. If the form
$L_\theta$ is definite, then $(M, \theta)$ is said to be a strictly
pseudoconvex pseudohermitian manifold. When $M$ is orientable, it
follows that the bundle of annihilators $H(M)^{\bot} := \{ \theta
\in T^*(M) :\; H(M) \subset \ker{\theta}\}$ admits a natural
orientation. We say that $\theta \in H(M)^{\bot}$ is positive, when
$L_\theta$ is positive definite. In this case, there is an
endomorphism $J:H(M) \rightarrow H(M)$ such that $J^2=-1$ and a unique
real vector field $T$, known as characteristic vector
field of $\theta$, such that $\theta(T)=1$ and $d\theta(T,\cdot)=0$. In particular,
it is possible to define a Riemannian metric $g_{\theta}$, known as
Webster metric, and a pseudohermitian connection $\nabla$ on
$TM \otimes \C$, known as Tanaka-Webster connection, compatible
with $J$ and $g_{\theta}$. From this connection, one gets the notions of
pseudohermitian curvature tensor, pseudohermitian Ricci tensor, torsion tensor on $T_{1,0}(M)$ and Webster scalar curvature of $\theta$. The Webster scalar curvature is the corresponding CR of the Riemannian
scalar curvature and was independently introduced by Webster \cite{Web} and Tanaka \cite{Tan}. Indeed, Webster constructed a
pseudohermitian curvature tensor similar to the Riemann curvature tensor, except that possibly admits torsion. Such a construction led him readily to the notions of pseudohermitian Ricci curvature tensor and to the so-called Webster scalar curvature. Throughout this paper, a CR manifold will mean an
orientable CR manifold endowed with a strictly pseudoconvex pseudohermitian structure.

A usual strategy in conformal geometry consists in choosing a metric in the conformal class of a fixed metric $g$ so that the geometry becomes simpler in some aspects. For instance, we have:

\begin{center}
\textbf{The Yamabe Problem.} Given a compact Riemannian manifold $(M,g)$ of dimension $n$, finding a metric conformal to $g$ with constant scalar curvature.
\end{center}

The great interest in this problem lies in the fact that its truth may mean that some topological questions can be reduced to other geometric ones on constant curvature models. For instance, when $n = 2$, the Yamabe problem is equivalent to the uniformization theorem from complex analysis. When $n \geq 3$, the Yamabe problem was completely solved by Aubin \cite{Au}, Schoen \cite{S} and Trudinger \cite{Tru} from the analytical formulation of a geometric invariant introduced by Yamabe \cite{Yam}, the so-called conformal Yamabe invariant.

A parallel problem was proposed in the CR manifold context by Jerison and Lee in \cite{JL2}. Precisely:

\begin{center}
\textbf{The CR Yamabe Problem.} Given a compact CR manifold $(M, \theta)$ of dimension $2n+1 \geq 3$, finding a pseudohermitian structure conformal to $\theta$ with positive orientation and constant Webster scalar curvature.
\end{center}

The CR Yamabe problem was partially solved by Jerison and Lee, within
four pioneer works \cite{JL1}, \cite{JL2}, \cite{JL3}, \cite{JL4},
for a compact CR manifold of dimension greater than $3$ and
non-locally CR-equivalent to the CR sphere. In short, given a
compact CR manifold $M$ of dimension $2n + 1$, they proved the
existence of a conformal geometric invariant $\lambda(M)$, analogous
to the conformal Yamabe invariant, satisfying:

\begin{enumerate}
\item[(a)] $\lambda(M)$ depends only on the CR structure on $M$;
\item[(b)] $\lambda(M) \leq \lambda(\s^{2n + 1})$, where $\s^{2n + 1}$ denotes the CR sphere in $\co^{n+1}$;
\item[(c)] if $M$ is non-locally CR-equivalent to $\s^{2n + 1}$ and $n > 1$, then $\lambda(M) < \lambda(\s^{2n + 1})$;
\item[(d)] if $\lambda(M) < \lambda(\s^{2n + 1})$, then $M$ admits a conformal pseudohermitian structure with constant Webster scalar curvature.
\end{enumerate}

\n The remaining cases, namely when $n = 1$ or the CR manifold is locally CR-equivalent to $\s^{2n + 1}$ for $n \geq 2$, were completed by Gamara and Yacoub in the works \cite{Ga} and \cite{GaY}.

Before we go further, a little bit of notation and background about Riemannian geometric tensors should be introduced. Let $(M,g)$ be a Riemannian manifold $(M,g)$ of dimension $n \geq 3$. Denote by $Rm_g$ the Riemann curvature tensor with respect to the metric $g$. It is well known that the tensor $Rm_g$ can be decomposed as

\[
Rm_g = W_g + S_g \odot g\, ,
\]

\n where $W_g$ stands for the Weyl conformal tensor and $S_g$ for the Schouten tensor

\[
S_g = \frac{1}{n-2} ( Ric_g - \frac{R}{2(n-1)} g )\, ,
\]

\n where $Ric_g$ and $R$ denote, respectively, the Ricci tensor and the scalar curvature of $g$ and $\odot$ represents the Kulkarni-Nomizu product (for instance, see \cite{Be}). This factorization plays a key role in Riemannian geometry because the Weyl tensor $W_g$ is a conformal invariant, being identically null if, and only if, either $n = 3$ or $n \geq 4$ and $M$ is locally conformally flat. In particular, the study of conformal deformations can be carried from the properties of the Schouten tensor $S_g$.

Let $\sigma_k$ be the $k$-th elementary symmetric function. For each symmetric real matrix $A$ of order $n$, denote the spectrum of $A$ by $(\lambda_1, \ldots, \lambda_n)$ and define

\[
\sigma_k(A) = \sum_{I_1 < \ldots < I_k} \lambda_{I_1} \ldots \lambda_{I_k}\, .
\]

\n Of course, $\sigma_k$ is invariant by the conjugation $A \mapsto P^t A P$, where $P$ is an orthogonal matrix, $\sigma_1(A) = {\rm trace}(A)$ and $\sigma_n(A) = \det(A)$.

In \cite{Via1}, Viaclovsky considered the function $\sigma_k$ on tensors of type $(1,1)$ and introduced the notion of $k$-curvature of a metric $g$ on $M$ as

\[
\sigma_k(g) := \sigma_k(g^{-1} S_g)\, ,
\]

\n where $g^{-1} S_g$ is locally seen as $(g^{-1} S_g)^i_j = g^{ik} (S_g)_{kj}$. For $k=1$, we have

\[
\sigma_1(g) = \frac{R}{2(n-1)}\, ,
\]

\n so that, modulo a positive constant, the $k$-curvature notion extends the scalar curvature one of $g$.

At the beginning of the 2000s, Viaclovsky \cite{Via1} and Chang, Gursky and Yang \cite{CGY} focused, independently, their attention on the following extension of the Yamabe problem:

\begin{center}
\textbf{The $k$-Yamabe Problem.} Given a compact Riemannian manifold $(M,g)$ of dimension $n$, finding a metric conformal to $g$ with constant $k$-curvature.
\end{center}

\n Since then, important results have been established in a series of outstanding works. We refer again to \cite{Au}, \cite{S}, \cite{Tru}, \cite{Yam} for the complete solution of the problem in the case $k=1$ from 1960 to 1984 and to \cite{BG}, \cite{CGY}, \cite{GeW}, \cite{GE}, \cite{GuW}, \cite{GV}, \cite{Li}, \cite{STW}, \cite{Via3}, \cite{Tru1} for the partial solution in the case $k > 1$ from 2000 to 2010.

In order to state two central results due to Viaclovsky, consider a
compact Riemannian manifold $(M,g)$ of dimension $n \geq 3$ and a
metric $\tilde{g}$ conformal to $g$, that is $\tilde{g} = e^{-2u} g$
for some smooth function $u$. After some computations, one arrives
at the following relation between the Schouten tensors
$S_{\tilde{g}}$ and $S_{g}$:

\[
S_{\tilde{g}} = S_{g} + \nabla^2 u + du \otimes du - \frac{|\nabla u|^2}{2} g\, ,
\]

\n so that $\sigma_k(\tilde{g})$ is constant if, and only if, $u$ satisfies the equation

\begin{equation} \label{FNE}
\sigma_k(S_{g} + \nabla^2 u + du \otimes du - \frac{|\nabla u|^2}{2} g) = \lambda e^{-2ku}\ \ \mbox{on}\ \ M
\end{equation}

\n for some real constant $\lambda$. It is worth mentioning that (\ref{FNE}) is a $k$-Hessian type fully nonlinear equation whenever $k > 1$ and that clearly becomes the Yamabe equation when $k=1$.

By (\ref{FNE}), one easily checks that the function $v$ defined by the relation $e^{-2u}=v^{p-2}$, with $p=2n/(n-2)$, satisfies the equation

\begin{equation} \label{FNE2}
\sigma_k ( {\dfrac{2}{(n-2)} \dfrac{\;V}{v}} ) = \lambda v^k\ \ \mbox{on}\ \ M\, ,
\end{equation}

\n where

\[
V = \dfrac{n-2}{2}\, v S_{g} -\nabla^2 v + \dfrac{n}{n-2} \dfrac{\nabla v \otimes \nabla v}{v} - \dfrac{1}{n-2} \dfrac{\|\nabla v\|^{2}_{g}}{v} g\, ,
\]

\n or equivalently,

\begin{equation} \label{FNE3}
{\cal L}_{k}[v]:=v^{{(1-k)}\frac{n+2}{n-2}}\sigma_k\lp{V}\rp = \lambda v^{\frac{n+2}{n-2}}\ \ \mbox{on}\ \ M\, .
\end{equation}

\n The operator ${\cal L}_k$ is naturally called the \textbf{$k$-Yamabe operator}. When $k=1$, the operator ${\cal L}_1$ is, up to a constant factor, equal to the Yamabe conformal operator and (\ref{FNE3}) yields the Yamabe equation. When $k > 1$, the $k$-Yamabe problem (or the equation (\ref{FNE3})) was first investigated by Viaclovsky in the handsome work \cite{Via1}. We highlight below some of his main results. The first one provides an answer about the variational nature of the equation (\ref{FNE}) and the second one concerns with a uniqueness result for the $k$-Yamabe problem.

Precisely, we have:

\vs{0.5cm}

\n {\bf Theorem A}\; (\textbf{Theorem 1 of \cite{Via1}}).
\emph{Let $(M,g)$ be a compact Riemannian manifold of dimension $n \geq 2$ and ${\cal M} = \{\tilde{g} \in [g]: \; \int_M dv_{\tilde{g}} = 1\}$. If $k \neq \frac{n}{2}$ and $(M,g)$ is locally conformally flat, then a metric $\tilde{g} \in {\cal M}$ is a critical point of the functional}

\[
{\cal F}_k : \tilde{g} \mapsto \int_M \sigma_k(\tilde{g})\;  dv_{\tilde{g}}
\]

\n \emph{constrained to ${\cal M}$ if, and only if, $\sigma_k(\tilde{g}) = \lambda_k$ for some constant $\lambda_k$. If $(M,g)$ is non-locally conformally flat, then the statement is true for $k = 1$ and $k = 2$.}

\vs{0.5cm}

\n {\bf Theorem B}\; (\textbf{Theorem 3 of \cite{Via1}}).
\emph{Let $(M,g)$ be a compact Riemannian manifold of dimension $n \geq 2$ with unit volume and non-null constant sectional curvature. Then, for any $k \in \{1, \ldots, n-1\}$, $g$ is the unique solution of unit volume of the equation $\sigma_k(\tilde{g}) = c$ in its class conformal unless $(M,g)$ is isometric to the round sphere $\s^n$. In this case, we have an $(n+1)$-parameter family of solutions that are the images of the standard metric under conformal diffeomorphisms of $\s^n$.}

\vs{0.5cm}

Consider the G\"{a}rding cone in $\R^n$

\[
\Gamma_k^+ = \{\Lambda = (\lambda_1, \ldots, \lambda_n) \in \R^n :\; \sigma_j(\Lambda) > 0\ \mbox{for all}\ j \leq k\}\, .
\]

\n A metric $g$ on $M$ is said to be $k$-positive if $\sigma_j(g)(x) > 0$ for all $x \in M$ and $1 \leq j \leq k$. In this case, we simply denote $g \in \Gamma_k^+$. If both $\tilde{g} = e^{-2u} g$ and $g$ are $k$-positive, we then say that $u$ is $k$-admissible.

In a precise way, we present another statement of the $k$-Yamabe problem.

\begin{center}
\textbf{A second version of the $k$-Yamabe Problem.} Given a compact Riemannian manifold $(M,g)$ of dimension $n$ with $k$-positive metric, finding a $k$-positive metric conformal to $g$ with constant $k$-curvature.
\end{center}

\n A fact that deserves mention, which has been recently proved, it is that certain Riemannian metrics $g$ on $M$ admit a
$k$-positive conformal metric, see \cite{GE}, \cite{Guan} and \cite{WS}. Thus, in some situations, the assumption of the
above problem can be removed.

Another result proved by Viaclovsky is that the equation (\ref{FNE}) is elliptic on any $k$-admissible solution $u$. Such an ellipticity played an important role in the work \cite{STW}.

There is a number of key results concerning with the fully nonlinear
equation (\ref{FNE}). We select some of them connected with the
existence of solution for the $k$-Yamabe problem. Namely, when $k=n$, under a few conditions, Viaclovsky
\cite{Via3} proved the existence of a solution for the problem
(\ref{FNE}). The interesting case $n=2k=4$ was studied by
Chang, Gursky and Yang \cite{CGY}. When the manifold is locally
conformally flat, the problem was independently solved by Guan and
Wang \cite{GuW} and Li and Li \cite{Li}. The case $k> n/2$ was
dealt by Gursky and Viaclovski \cite{GV}. When $k=2$ and the
manifold is non-locally conformally flat, the problem was solved by
Ge and Wang \cite{GeW}. An interesting strategy in solving the
remaining cases, namely $2\leq k \leq n/2$ and non-locally conformally flat manifolds, was presented by Sheng, Trudinger and
Wang in \cite{STW}, whose solution was given under the assumption of
the problem (\ref{FNE}) be variational. Unfortunately, just in these
cases, the problem is not variational as proved by Branson and Gover
in \cite{BG}, leaving the question still open.

In the present paper, we introduce the $k$-Yamabe problem in the CR manifolds context and discuss on some related issues. For instance, we establish the corresponding of Theorem A in this new setting and, from the variational property, we produce a partial version of Theorem B on the CR sphere.

Let $(M, \theta)$ be a CR manifold of dimension $2n+1$. The Chern pseudo-conformal tensor \cite{CM} is a well-known CR geometric invariant which is null if, and only if, either $n=1$ or $M$ is locally CR-equivalent to the CR sphere. The similarity between this tensor and the Weyl conformal tensor leads naturally to ask about the CR analogue of the Schouten tensor. Inspired on the work of Chern and Moser \cite{CM}, Webster \cite{Web} noticed that the pseudohermitian curvature tensor $R_{\theta}$, seen as a tensor of type $(4,0)$, can be decomposed as

\begin{equation}
R_{\theta} = Ch_{\theta} + S_{\theta} \boxdot L_{\theta}\, ,
\end{equation}

\n where $Ch_{\theta}$ stands for the Chern pseudo-conformal tensor, $\boxdot$ is the CR corresponding of the Kulkarni-Nomizu product and $S_{\theta}$ is the tensor of type $(2,0)$ given by

\begin{equation}
S_{\theta} = \frac{1}{n+2} ( Ric_{\theta} - \frac{R}{2(n+1)} L_{\theta} )\, ,
\end{equation}

\n where $Ric_{\theta}$ and $R$ denote, respectively, the pseudohermitian Ricci tensor and the Webster scalar curvature. The tensor $S_{\theta}$ will be called the pseudohermitian Schouten tensor associated to the pseudohermitian structure $\theta$. In parallel to what occurs in the Riemannian geometry, the pseudohermitian Schouten tensor plays an important role in the CR geometry.

Since the tensor $S_{\theta}$ is hermitian on $T_{1,0}(M)$, we can introduce the notion of pseudohermitian $k$-curvature $\sigma_k(\theta)$ of $\theta$ as

\[
\sigma_k(\theta) := \sigma_k({L_{\theta}}^{-1} S_\theta)\, ,
\]

\n where ${L_{\theta}}^{-1} S_{\theta}$ is locally given by $S_{\alpha}\,^{\beta}:=h^{\beta \oga} (S_\theta)_{\alpha \oga}$, being $h^{\beta \oga}$ and $(S_\theta)_{\alpha \oga}$, respectively, the components of the tensors ${L_{\theta}}^{-1}$ and $S_\theta$. As we shall see,

\[
\sigma_1(\theta) = \frac{R}{2(n+1)}\, .
\]

\n In other words, the pseudohermitian $1$-curvature is a constant multiple of the Webster scalar curvature.

In a natural way, we propose the following question:

\begin{center}
\textbf{The CR $k$-Yamabe Problem.} Given a compact CR manifold $(M, \theta)$ of dimension $2n+1 \geq 3$, finding an orientation preserving pseudohermitian structure conformal to $\theta$ with constant pseudohermitian $k$-curvature.
\end{center}

An orientation preserving pseudohermitian structure $\tilde{\theta}$ conformal to $\theta$ can be written as $\tilde \theta = e^{2u} \theta$ for some smooth function $u$ on $M$. The relation between the pseudohermitian Schouten tensors $S_{\tilde \theta}$ and $S_\theta$ is given in terms of $u$ by

\begin{equation} \label{1}
S_{\tilde \theta} = S_\theta -2 \, {\nabla}^{2}u+{\lp {iT_{\theta} u-{\normc {du} 2 \theta} }\rp} L_{\theta} \, ,
\end{equation}

\n where $\nabla^2 u$ denotes the complex Hessian of $u$ and $T_{\theta}$ stands for the characteristic direction associated to $d \theta$. From (\ref{1}), it follows that $\sigma_k(\tilde{\theta})$ is constant if, and only if, $u$ satisfies

\begin{equation} \label{FNECR}
\sigma_k(S_\theta -2 \, \nabla^2 u + {\lp {iT_{\theta} u-{\normc {du} 2 \theta} }\rp}L_{\theta}) = \lambda e^{2ku}\ \ \mbox{on}\ \ M
\end{equation}

\n for some constant $\lambda$. Remark that (\ref{FNECR}) is a $k$-Hessian type fully nonlinear equation whenever $k > 1$ and that recovers the CR Yamabe equation when $k = 1$.

Performing now the change $e^{2u} = v^{p-2}$, with $p=2+2/n$, one obtains

\begin{equation} \label{FNECR2}
\sigma_k ({\dfrac{2}{n} \dfrac{\;V}{v}}) = \lambda v^k\ \ \mbox{on}\ \ M\, ,
\end{equation}

\n where

\[
V = \dfrac{n}{2} \, v S_{\theta} -\nabla^2 v + \dfrac{n}{n-2} \dfrac{\nabla v \otimes \nabla v}{v} +\dfrac{1}{2} ( {iT_{\theta} v-\dfrac{\|\nabla v\|^{2}_{\theta}}{nv}} ) L_{\theta}\, .
\]

\n Moreover, $v$ also satisfies

\begin{equation} \label{FNECR3}
{\cal L}_{k}[v]:=v^{{(1-k)}\frac{n+2}{n}}\sigma_k ( {V} ) = \lambda v^{\frac{n+2}{n}}=\lambda v^{p-1}\ \ \mbox{on}\ \ M\, .
\end{equation}

\n The operator ${\cal L}_k$ will be naturally called the \textbf{$k$-Yamabe operator}. As we shall see in the next section, the operator ${\cal L}_1$ is a multiple of the CR Yamabe conformal operator introduced by Jerison and Lee \cite{JL2}.

Regarding again the G\"{a}rding cone $\Gamma_{k}^{+}$, we say that a
pseudohermitian structure $\theta$ on $M$ is $k$-positive if
$\sigma_{j}(\theta)(x)>0$ for all $x \in M$ and $1 \leq j \leq k$.
When both $\tilde{\theta} = e^{2u} \theta$ and $\theta$ are
$k$-positive, we say that $u$ is $k$-admissible. Denote by
$[\theta]_{+}$ the set of the $k$-positive pseudohermitian
structures compatible with $\theta$.

We now state a variation of the CR $k$-Yamabe problem.

\begin{center}
\n \textbf{A second version of the CR $k$-Yamabe problem.} Given a compact CR manifold $(M, \theta)$
of dimension $2n+1 \geq 3$ with $k$-positive pseudohermitian structure, finding a $k$-positive pseudohermitian
structure conformal to $\theta$ with constant pseudohermitian $k$-curvature.
\end{center}

Motivated by previous contributions on the Riemannian $k$-Yamabe problem, a strategy in solving the above problem is knowing if the problem is variational. In order this, we consider the following functional

\begin{equation} \label{2}
Y_{k}: \tilde \ta \mapsto \int_{M} \si_k(\tilde \ta) \, dV_{\tilde \ta}\, ,
\end{equation}

\n constrained to the subspace

\[
{\mathcal M}=\{\tta \in [\ta]: \int_M dV_{\tta} =1\}
\]

\n of conformal pseudohermitian structures $\tta$ of unit volume and compatible with $\theta$. From the relation between ${\cal L}_1$ and the CR Yamabe conformal operator, it follows that $Y_1(\tilde \ta)=\frac{1}{2(n+1)}Y(\tilde \ta)$, where $Y$ is the previously mentioned CR Yamabe functional associated to $\ta$.

Before presenting our main results, we need to recall some geometric concepts in CR geometry.

Given a CR manifold $(M, \ta)$ of dimension $2n + 1$, the pseudohermitian Cotton tensor $C_\ta$ is a tensor of type $(3,0)$ locally defined by

\begin{equation}
C_{\rro \osi \ga} = S_{\rro \osi ; \ga} - S_{\ga \osi ; \rro}
\end{equation}

\n where $S_{\rro \osi}$ represents the components of the pseudohermitian Schouten tensor $S_{\ta}$.

Our first result concerns to a variational property of the problem (\ref{FNECR}).

\begin{teorema}
Let $(M, \ta)$ be a compact CR manifold of dimension $2n+1$. If $k \geq 1$, then a conformal pseudohermitian structure $\tilde{\ta} \in {\cal M}$ satisfying $C_{\tilde \ta} = 0$ is a critical point of (\ref{2}) if, and only if, $\sigma_k (\tilde \ta) = \lambda_k$ for some constant $\lambda_k$.
\end{teorema}

The above theorem will be a fundamental tool in the study of the CR
$k$-Yamabe problem on the CR sphere, namely in the proof of Corollary 1.4
and Theorem 1.4. Indeed, the example of the sphere is the what
inspired us in the search for a variational property as above. On
the other hand, it would be very interesting the construction of
other examples on which the above theorem can be applied.

We now recall in brief the notion of pseudohermitian torsion tensor. Consider $n$ forms $\{\theta^1, \ldots, \theta^n\}$ of type $(1,0)$ whose restrictions to $T_{1,0}(M)$ form a basis for $T_{1,0}^{*}(M)$ and such that $T_\theta \,\rfloor \, \theta^\alpha = 0$ for
$\alpha=1,\ldots,n$. In other words, $\{\theta^{\alpha}\}$ is an admissible coframe. Note that $\{\ta, \ta^1, \ldots, \ta^n, \ta^{\overline{1}},
\ldots, \ta^{\overline{n}} \}$ is a coframe for $T(M) \otimes \C$. Choose an admissible coframe $\{\ta^\al\}$ and consider its dual frame $\{T_{\al}\}$ for $T_{1,0}(M)$ with the convention $\overline{T_\al}=T_{\oal}$. Webster showed that there are uniquely determined $1$-forms $\om_{\ga}^{\;\; \al}$ and $\tau^{\al}$ on $M$ satisfying
\[
d\ta^{\al}=\ta^{\ga} \wedge \om_{\ga}^{\  \al} + \ta \wedge \tau^{\al} \, ,
\]

\[
dh_{\al \obe}=\om_{\al \obe}+\om_{\obe \al} \, ,
\]

\begin{equation} \label{simtor}
\tau_{\al} \wedge \ta^\al = 0 \, ,
\end{equation}

\n where $h_{\al \obe}$ demote the components of the Levi form $L_\theta$. By (\ref{simtor}), we can write

\[
    \tau_{\al}=A_{\al \ga}\ta^\ga
\]

\n with $A_{\al \ga}=A_{\ga \al}$.

The tensor whose components are $A_{\al \be}$ is called the pseudohermitian torsion of $\theta$ and denoted by $\tau_\theta$. The pseudohermitian torsion is said to be parallel with respect to the Tanaka-Webster connection $\nabla$, if $\nabla \tau_\theta = 0$. An immediate consequence of the
Bianchi identities (see \cite{Lee2}) is that parallel torsion implies in vanishing of the Cotton tensor.

\begin{coro}
Let $(M, \ta)$ be a compact CR manifold of dimension $2n+1$. If $k \geq 1$, then a conformal pseudohermitian structure $\tilde{\ta} \in {\cal M}$ with parallel torsion is a critical point of (\ref{2}) if, and only if, $\sigma_k (\tilde \ta) = \lambda_k$ for some constant $\lambda_k$.
\end{coro}

Another particularly interesting situation holds when the CR manifold $(M, \theta)$ is torsion-free. In this case, if $\tilde \ta =  e^{2u} \ta$ is a pseudohermitian structure conformal to $\theta$, under the Cotton-admissibility condition for $u$, we have $C_{\tilde \ta} = 0$. Indeed, a function $u$ is said to be Cotton-admissible if satisfies

\[
u_{\al} u_{\obe \si} - u_{\si} u_{\obe \al} = 0
\]

\n for all $\al , \be , \si \in \{1, ... ,n\}$. A standard class of Cotton-admissible functions is given by CR-pluriharmonic functions,
precisely functions that are real part of CR functions, see \cite{Lee2}.

Consider the sets

\[
\mathcal{C}[\ta] = \left\{e^{2u}\ta:u \; \mbox{is Cotton-admissible}\right\}\, ,
\]

\[
\mathcal{C}[\ta]_{+} = \left\{e^{2u}\ta:u \; \mbox{is $k$-admissible and Cotton-admissible}\right\}\, .
\]

\begin{coro}
Let $(M, \ta)$ be a torsion-free compact CR manifold of dimension $2n+1$. If $k \geq 1$, then a conformal pseudohermitian structure $\tilde{\ta} \in {\cal M} \cap \mathcal{C}[\ta]$ is a critical point of (\ref{2}) if, and only if, $\sigma_k (\tilde \ta) = \lambda_k$ for some constant $\lambda_k$.
\end{coro}

The next result requires the notion of pseudo-Einstein structure. We recall that a pseudohermitian structure $\theta$ on $M$ is said to be pseudo-Einstein if, constrained to $T_{1,0}(M)$, one has $Ric_{\ta}=\frac{R}{n} L_{\ta}$.

\begin{teorema}
Let $(M, \ta)$ be a compact CR manifold of dimension $2n+1$. Assume that $k \geq 1$ and $\theta$ is pseudo-Einstein. If $\tilde \ta \in [\ta]_{+}$ has constant pseudohermitian $k$-curvature and null Cotton tensor, then $\tilde \ta$ is pseudo-Einstein too.
\end{teorema}

An immediate consequence of the preceding result is

\begin{coro}
Let $(M, \ta)$ a compact CR manifold of dimension $2n+1$. Assume that $k \geq 1$ and $\theta$ is torsion-free and pseudo-Einstein. If $\tilde \ta \in \mathcal{C}[\ta]_{+}$ has constant pseudohermitian $k$-curvature, then $\tilde \ta$ is pseudo-Einstein too.
\end{coro}

We also obtain a classification result of solutions for the $k$-Yamabe problem on the CR sphere $\s^{2n + 1}$, the so-called Obata type result. In fact, we have:

\begin{coro}
Let $(\s^{2n+1}, \hat{\ta})$ be the CR sphere of dimension $2n + 1$. Let $k \geq 1$. If the pseudohermitian structure $\tilde \ta = e^{2u} \hat{\ta} \in \mathcal{C}[\hat{\ta}]_{+}$ has constant pseudohermitian $k$-curvature, then $\tilde \ta$ is equal to a multiple of the canonical form $\hat{\ta}$ by some CR-automorphism of $\s^{2n+1}$.
\end{coro}

In view of well-known results on the Riemannian $k$-Yamabe problem and the CR Yamabe problem, we believe that Corollary 1.4 can be extended as follows:

\begin{conj}
Let $(\s^{2n+1}, \hat{\ta})$ be the CR sphere of dimension $2n + 1$. Let $k \geq 1$. If the pseudohermitian structure $\tilde \ta = e^{2u} \hat{\ta} \in [\hat{\ta}]_{+}$ has constant pseudohermitian $k$-curvature, then $\tilde \ta$ is equal to a multiple of the canonical form $\hat{\ta}$ by some CR-automorphism of $\s^{2n+1}$.
\end{conj}

\n Of course, the truth of Conjecture 1.1 would imply in the complete solution of the $k$-Yamabe problem on $\s^{2n+1}$.\\

Given a CR manifold $(M, \ta)$ of dimension $2n + 1$, we introduce the following constants:

\[
\la_k(M) = \inf \{  Y_k(\tta) : \tta \in {\cal M} \cap [\ta]\}\, ,
\]

\[
\la_k^{+}(M) = \inf \{  Y_k(\tta) : \tta \in {\cal M} \cap [\ta]_{+}\}
\]

\n and

\[
\la_k^{\mathcal{C}}(M) = \inf \{  Y_k(\tta) : \tta \in {\cal M} \cap \mathcal{C}[\ta]_{+}\}\, .
\]

\n In the case that $M$ is compact, each above constant is finite and, moreover,

\[
\la_k(M) \leq \la_k^{+}(M) \leq \la_k^{\mathcal{C}}(M)\, .
\]

\n These constants will be referred as the CR $k$-Yamabe constants. A natural issue is knowing if, at least, one of the above inequalities is strict for some compact CR manifold. This is an open question.

An important consequence of Theorem 1.1 and Corollary 1.4 is that the above introduced constants, in the case of the CR sphere $(\s^{2n+1}, \hat{\ta})$, are attained by a multiple of the form $\hat{\ta}$ by some CR-automorphism of $\s^{2n+1}$. Using this fact, we produce the following result:

\begin{teorema}
Let $(\s^{2n+1}, \hat{\ta})$ be the CR sphere of dimension $2n + 1$. Then, for any $k \geq 1$,

\[
\la_k(\s^{2n+1}) = \la_k^{+}(\s^{2n+1}) = \la_k^{\mathcal{C}}(\s^{2n+1}) = C(n,k)\pi^{k}\, ,
\]

\n where $C(n,k)$ denotes the Newton binomial coefficient given by $n!/k!(n-k)!$.
\end{teorema}

Another interesting application of Theorem 1.1 and Corollary 1.4 is the existence of extremal functions for Hessian type Folland-Stein inequalities on the Heisenberg group. We refer to \cite{BCM2} for more details.

In parallel to what occurs in the Riemannian context, we also hope that the CR $k$-Yamabe constants play a central role on the existence of solution for the CR $k$-Yamabe problem. Our next theorem establishes an upper bound for the CR $k$-Yamabe constant $\la_k^{+}(M)$ on any compact CR manifold $(M, \theta)$ of dimension $2n + 1$ and also a classification result in the equality case.

Precisely, we have:

\begin{teorema}
Let $(M, \ta)$ be a compact CR manifold of dimension $2n+1$. Then, for any $k \geq 1$,

\[
\la_{k}^{+}(M) \leq \la_{k}^{+}(\s^{2n+1})\, .
\]

\n Furthermore, assuming $\lambda(M) \geq 0$, the above inequality yields equality if, and only if, $M$ is locally CR-equivalent to the CR sphere $\s^{2n+1}$.
\end{teorema}

Thanks to the previous result, to the solution of the Riemannian and CR Yamabe problems and, specially, to developments due to Sheng, Trudinger and Wang \cite{STW} on the $k$-Yamabe problem and to the recent work due to Li and Monticelli \cite{LiMo} on fully nonlinear PDEs in the Heisenberg group, we are led to raise the following conjecture:

\begin{conj}
Let $(M, \ta)$ be a compact CR manifold of dimension $2n+1$. If $k \geq 1$ and

\[
\la_{k}^{+}(M) < \la_{k}^{+}(\s^{2n+1})\, ,
\]

\n then the CR $k$-Yamabe constant $\la_{k}^{+}(M)$ is attained by a positive smooth function $u$. In particular, the CR $k$-Yamabe problem has a solution provided that it is variational and $\la_{k}^{+}(M) < \la_{k}^{+}(\s^{2n+1})$.
\end{conj}

\n Indeed, we believe that this conjecture can be partially achieved with aid of the program introduced in \cite{STW} and of some ideas developed in the work \cite{LiMo}.

The proof of Theorem 1.4 relies on an ellipticity result for the Hessian type fully nonlinear equation (\ref{FNECR}). In a precise way, we have:

\begin{teorema}
Let $(M, \ta)$ be a compact CR manifold with $k$-positive pseudohermitian structure $\theta$ and be $u \in C^\infty(M)$ a positive function such that

\[
{\cal L} u := \si_k(S_{\ta} -2\grad^{2}u+(iT_{\ta} u-\normc {du} 2 \ta )L_{\ta}) > 0 \ \ \mbox{on}\ \ M\, .
\]

\n Then ${\cal L}$ is elliptic at $u$, that is the linearized of ${\cal L}$ at $u$ provides an elliptic operator.
\end{teorema}

The organization of paper is as follows. In Section 2, we present some geometric tensors which are essential in CR geometry and the notion of pseudohermitian $k$-curvature. Section 3 is devoted to the formulation of the $k$-Yamabe problem on compact CR manifolds. In Section 4, we provide the complete proof of Theorems 1.1 and 1.2. In Section 5, we focus on the proof of Theorem 1.5. Finally, in Section 6, we present the proof of Theorems 1.3 and 1.4.\\

\section{Pseudohermitian geometric tensors and the CR $k$-curvature}

In this section, we recall some important geometric tensors in CR
geometry and introduce the $k$-curvature notion. The strategy
of construction of this curvature is inspired in the work
\cite{Via1} of Viaclovsky.

Our first goal is to search for a tensor associated to CR manifolds that
plays the same role exerted by the Schouten tensor in Riemannian geometry. Such a CR tensor will be naturally called pseudohermitian Schouten tensor.

Unless otherwise said, we adopt the following conventions for the range of indices in this work. The Greek indices $\alpha, \beta, \gamma, \ldots$ vary from $1$ to $n$, while the Latin indices $A, B, C, \ldots$ vary in $\{0, 1, \ldots , n, \overline{1}, \ldots, \overline{n}\}$, with the
convention $T_0 = T$. We also adopt the Einstein summation convention and use the matrix $h_{\alpha \overline{\beta}}$, formed by the components of the Levi form $L_\theta$, to raise and lower indices.

Let $(M,\theta)$ be a CR manifold of dimension $2n + 1$. The
\textbf{pseudo-conformal curvature tensor} (or \textbf{Chern
tensor}), introduced by Chern and Moser in \cite{CM}, is the tensor
$Ch_\theta$ of type $(3,1)$ with components

\begin{equation} \label{Chern}
C_{\be \;\; \la \osi}^{\;\; \al} = R_{\be \;\; \la \osi}^{\;\; \al} - \frac{1}{n+2} ( {R_{\be}^{\;\; \al}h_{\la \osi}
+ R_{\la}^{\;\; \al}h_{\be \osi} +\dl^{\al}_{\be}R_{\la \osi}+\dl^{\al}_{\la}R_{\be \osi}} ) + \frac{R}{(n+1)(n+2)} ( {
\dl^{\al}_{\be}h_{\la \osi}+\dl^{\al}_{\la}h_{\be \osi} } )\, ,
\end{equation}

\noindent where $R_{\be \;\; \la \osi}^{\;\; \al}$, $R_{\be \omi}$, $h_{\al \obe}$ and $R$ denote, respectively, the components of the
pseudohermitian curvature tensor, of the pseudohermitian Ricci tensor, of the Levi form and the Webster scalar curvature associated to the pseudohermitian structure $\ta$. Note that $C_{\be \;\; \la \osi}^{\;\; \al}$ are similar to the
components of the Weyl tensor, but viewed as a tensor of type $(3,1)$.

The Chern tensor $Ch_\theta$ is similar to the Weyl tensor of the Riemannian geometry in the sense that $Ch_\theta$ vanishes if, and only if, $n=1$ or $M$ is locally CR-equivalent to the CR sphere. Thus, it is natural to seek an analogous of the Schouten tensor starting from a decomposition of the pseudohermitian curvature tensor into two parts, being one of them given by the Chern tensor.

Performing a contraction in (\ref{Chern}) with the Levi form $h_{\al \omi}$, we find

\begin{eqnarray*}
C_{\be \omi \la \osi} = C_{\be \;\; \la \osi}^{\;\; \al}h_{\al \omi} &=& R_{\be \omi \la \osi}- \frac{1}{n+2} ( {R_{\be \omi}h_{\la \osi}  + R_{\la \omi}h_{\be \osi} + R_{\la \osi}h_{\be \omi}+R_{\be \osi}h_{\la \omi}} ) \\
 &&+ \frac{R}{(n+1)(n+2)} ( { h_{\la \osi}h_{\be \omi}+h_{\be \osi}h_{\la \omi} } ) \, .
\end{eqnarray*}

\noindent Consider the following product $\boxdot$ on
$(2,0)$-tensors associated to $(M, \theta)$:

\[
(K \boxdot S)(X,Y,Z,W) = K(X,Y)S(Z,W)+S(X,Y)K(Z,W) +
K(X,W)S(Z,Y)+S(X,W)K(Z,Y)\, .
\]

\noindent This product plays in CR geometry the same role played by the
Kulkarni-Nomizu product in Riemannian geometry.

Given a frame $\{T_{\al}\}$ associated to $(M, \ta)$, one has

\[
(Ric_{\ta} \boxdot L_{\ta})(T_{\be},T_{\omi},T_{\la},T_{\osi}) =
R_{\be \omi}h_{\la \osi}+R_{\la \omi}h_{\be \osi}+R_{\la \osi}h_{\be
\omi}+R_{\be \osi}h_{\la \omi} \, ,
\]

\[
(L_{\ta} \boxdot L_{\ta})(T_{\be},T_{\omi},T_{\la},T_{\osi}) = 2 (
h_{\la \osi}h_{\be \omi}+h_{\be \osi}h_{\la \omi} )\, .
\]

\noindent Then, we can write

\[
C_{\be \omi \la \osi}=R_{\be \omi \la \osi}- {\frac{1}{n+2} ( {Ric_{\ta}-\frac{R}{2(n+1)}L_{\ta}} ) }
\boxdot L_{\ta}(T_{\be},T_{\omi},T_{\la},T_{\osi})\, .
\]

\noindent Therefore, the \textbf{pseudohermitian Schouten tensor} is naturally defined as

\[
S_{\ta} = \frac{1}{n+2} ( {Ric_{\ta}-\frac{R}{2(n+1)}L_{\ta}} )\, ,
\]

\noindent since

\[
R_{\theta} = Ch_\theta + S_{\ta} \boxdot L_{\ta}\, .
\]

\noindent Note that $S_\theta$ is a $(2,0)$-tensor on $(M, \theta)$ whose components are given by

\begin{equation} \label{schouten}
S_{\al \obe} = \frac{1}{n+2} ( {R_{\al \obe}- \frac{R}{2(n+1)} h_{\al \obe}} )\, .
\end{equation}

\noindent It's easy to see that $S_{\al \obe}$ is hermitian on $T_{1,0}(M)$. Using the natural complex extensions of $R_{\al \obe}$ and
$h_{\al \obe}$, we can consider $S_{\al \obe}$ defined on $T(M) \otimes \C$. Assuming this definition, $S_{\ta}$ inherits the properties

\[
S_{\al \obe}=S_{\obe \al } \;\;\; , \;\;\; \overline{S_{\al
\obe}}=S_{\oal \be}\, .
\]

\noindent In addition, the tensor $S_{\ta}$ has the following matrix representation with respect to a fixed frame:

\[
S_{\ta}:
\begin{bmatrix}
i\lp \frac{n-1}{n+2} \rp A_{\al \be} & S_{\oal \be} & \frac{1}{n+2}A_{\al \be ;}^{\;\;\;\;\;\; \si}\\
S_{\al \obe} & -i\lp \frac{n-1}{n+2} \rp A_{\oal \obe}  & \frac{1}{n+2}A_{\al \be ;}^{\;\;\;\;\;\; \osi}\\
0 & 0 & 0
\end{bmatrix}\, .
\]

\noindent where $A_{\al \be}$ denotes the components of the pseudohermitian torsion tensor as in the introduction. In particular, $S_{\al \obe}$ determines $S_{\ta}$ only when $(M, \theta)$ is torsion-free.

Since the Chern tensor is a conformal CR invariant, in order to investigate conformal deformations of pseudohermitian structures, we
focus our attention on the tensor with components $S_{\al \obe}$.

Let $\tta = e^{2u}\ta$ be a pseudohermitian structure conformal to $\ta$. Consider an admissible coframe $\{\ta^{\al}\}$ associated to $(M,\ta)$. For most purposes, it is appropriate to work with the coframe $\{\tta^{\al}=\ta^{\al}+2iu^{\al}\ta\}$, which is admissible for
$(M, \tta)$ and is dual to the original holomorphic frame. Regarding the coframe $\{\tta^\al\}$, the components of the Levi form $L_{\tta}$ satisfy ${\tilde h}_{\al \obe}=e^{2u} h_{\al \obe}$, see \cite{Lee2}. Furthermore,

\[
{\til R}_{\al \obe} = R_{\al \obe}-(n+2)(u_{\al \obe}+u_{\obe
\al})-(u^{\;\; \ga}_{\ga} + u^{\;\;
\oga}_{\oga}+4(n+1)u_{\ga}u^{\ga})h_{\al \obe}\, ,
\]

\[
{\til R} = e^{-2u}(R-2(n+1)(u^{\;\; \ga}_{\ga} + u^{\;\;
\oga}_{\oga})-4n(n+1)u_{\ga}u^{\ga})\, .
\]

\noindent Thus,

\begin{eqnarray*}
\til{S}_{\al \obe}&=&\frac{1}{n+2}({\til{R}_{\al \obe}-\frac{\til{R}}{2(n+1)}\tih_{\al \obe}})\\
&=&\frac{1}{n+2}\lp{R_{\al \obe}-(n+2)(u_{\al \obe}+u_{\obe \al})-(u^{\;\; \ga}_{\ga} + u^{\;\; \oga}_{\oga}+4(n+1)u_{\ga}u^{\ga})h_{\al \obe}}\rp\\
& &-\frac{1}{n+2} {\frac{e^{-2u}(R-2(n+1)(u^{\;\; \ga}_{\ga} + u^{\;\; \oga}_{\oga})-4n(n+1)u_{\ga}u^{\ga})}{2(n+1)}} e^{2u}h_{\al \obe}\\
&=&\frac{1}{n+2}R_{\al \obe}-(u_{\al \obe}+u_{\obe \al}) -\frac{u^{\;\; \ga}_{\ga} + u^{\;\; \oga}_{\oga}+4(n+1)u_{\ga}u^{\ga}}{n+2}h_{\al \obe}\\
& &-\frac{R}{2(n+1)(n+2)}h_{\al \obe}+\frac{u^{\;\; \ga}_{\ga} + u^{\;\; \oga}_{\oga}+2nu_{\ga}u^{\ga}}{n+2} h_{\al \obe}\\
&=& \frac{1}{n+2} ({R_{\al \obe} -\frac{R}{2(n+1)}h_{\al \obe}
})-(u_{\al \obe}+u_{\obe \al}) - 2 u_{\ga}u^{\ga} h_{\al \obe}\, ,
\end{eqnarray*}

\noindent so that

\[
\til{S}_{\al \obe}=S_{\al \obe}-(u_{\al \obe}+u_{\obe
\al})-2u_{\ga}u^{\ga}h_{\al \obe}\, .
\]

\noindent We now introduce a suitable hermitian form $L_{\ta}^{*}$
on $T^*(M) \otimes \C$. Let $\{\ta, \ta^{\al}, \ta^{\oal}\}$ be a
coframe on $T^*(M) \otimes \C$. Define

\[
L_{\ta}^{*}(\ta^{A} , \ta^{B}):=h^{AB}
\]

\noindent and extend linearly to all $T^*(M) \otimes \C$. In this
way,

\begin{eqnarray*}
&& L_{\ta}^{*}(\ta^{\al}, \ta^{\obe})=h^{\al \obe},\\
&& L_{\ta}^{*}(\ta^{\oal}, \ta^{\be})=h^{\oal \be}=\overline {h^{\al \obe}}=h^{\be \oal},\\
&& L_{\ta}^{*}(\ta^{\al},
\ta^{\be})=L_{\ta}^{*}(\ta^{\oal},\ta^{\obe}) = L_{\ta}^{*}(\ta,
\ta^{B})=L_{\ta}^{*}(\ta^{A}, \ta)=0\, .
\end{eqnarray*}

\noindent So, for any smooth function $u$ on $M$, we have

\begin{eqnarray*}
\normc {du} 2 \ta &=& L_{\ta}^{*}(du,du)\\
 &=&L_{\ta}^{*}(u_{\al}\ta^{\al} + u_{\oal}\ta^{\oal} +u_0 \ta , u_{\be}\ta^{\be} + u_{\obe}\ta^{\obe} +u_0 \ta)\\
 &=& u_{\al}u_{\obe}h^{\al \obe}+u_{\oal}u_{\be}h^{\oal \be}\\
&=&2u_{\al}u_{\obe}h^{\al \obe}\\
 &=& 2u_{\al}u^{\al}\, .
\end{eqnarray*}

\noindent Using the identity $u_{\al \obe} - u_{\obe \al} = i u_0
h_{\al \obe}$ (see \cite{Lee1}), we can state the following
proposition:

\begin{propo}
\label{Shouten1} Let $(M, \theta)$ be a CR manifold of dimension $2n + 1$ and $\tta = e^{2u}\ta$ be a pseudohermitian structure conformal
to $\theta$. Regarding the admissible coframe $\{\tta^{\al} = \ta^{\al}+2iu^{\al}\ta\}$ of $(M,\tta)$, we have

\begin{eqnarray*}
&& \til{S}_{\al \obe} = S_{\al \obe}-u_{\al \obe}-u_{\obe \al}-\normc {du} 2 \ta h_{\al \obe}\, ,\\
&& \til{S}_{\al \obe} = S_{\al \obe}-2u_{\al \obe}+(iu_0-\normc {du}
2 \ta) h_{\al \obe}\, .
\end{eqnarray*}

\n In particular, we have on $T_{1,0}(M)$,

\begin{equation} \label{schout}
S_{\tta} = S_{\ta} -2\grad^{2}u+(i T_\ta u-\normc {du} 2 \ta )L_{\ta}\, .
\end{equation}
\end{propo}

A tensor $K: T(M) \otimes \C \times T(M) \otimes \C \rightarrow \C$ of type $(2,0)$ is said to be pseudohermitian if, for any frame
$\{T_{\al}\}$, satisfies

\begin{enumerate}
    \item $K_{AB}=K_{\bar A \bar B} \, , \;\;\;\;$ ($K$ is real)
    \item $\overline {K_{\al \obe}}=K_{\be \oal}\, , \;\;\;\;$ ($K$ is Hermitian on $T_{1,0}$)
    \item $K_{AB}=K_{BA} \, , \forall A,B \in \{1,...,n,\bar{1}, ... , \bar{n}\} \;\;\;\;$ ($K$ is symmetric on $H(M) \otimes \C$)
\end{enumerate}

\n where $K_{AB}=K(T_A,T_B)$. The Levi form associated to a pseudohermitian structure is an example of pseudohermitian tensor with $h_{\al \be}=0$. Other ones are the pseudohermitian Schouten tensor $S_{\ta}$ and the pseudohermitian Ricci tensor $Ric_{\ta}$.

Since $L_{\ta}$ is nondegenerate on $H(M) \otimes \C$, for each $X
\in H(M) \otimes \C$, there exists a unique $X^* \in H(M) \otimes
\C$ such that

\[
K(X,\overline Y) = L_{\ta}(X^*,\overline Y)
\]

\noindent for all $Y \in H(M) \otimes \C$. Moreover, the map $X
\rightarrow X^*$ is $\C$-linear. Therefore, related to $K$, one
defines a tensor $K^*$ of type $(1,1)$ on $T(M) \otimes \C$ by

\[
K^* (X) = \left\{
\begin{array}{ll}
X^* & {\rm if}\ X \in H(M) \otimes \C\\
0 & {\rm otherwise}
\end{array}
\right.\; .
\]

\n The tensor $K^*$ is called \textbf{adjoint representation} of $K$
with respect to the Levi form $L_{\ta}$.

Taking a frame $\{T_{\al}\}$, we have

\begin{eqnarray*}
K_{\al \obe}&=&K(T_{\al}, T_{\obe})=L_{\ta}(K^* (T_{\al}), T_{\obe})\, ,\\
K_{\al \be}&=&K(T_{\al}, T_{\be})=L_{\ta}(K^* (T_{\al}), T_{\be})\,
.
\end{eqnarray*}

\noindent Writing

\[
K^*(T_{\al})=(K^*)_{\al}^{\;\; A}T_A\, ,
\]

\noindent one obtains

\begin{eqnarray*}
(K^*)_{\al}^{\;\; A}h_{A \obe} = K_{\al \obe} \; ,
\;\;(K^*)_{\al}^{\;\; A}h_{A \be} = K_{\al \be}\, .
\end{eqnarray*}

\noindent Consequently,

\begin{eqnarray*}
(K^*)_{\al}^{\;\; \ga}h_{\ga \obe} = K_{\al \obe} \; , \;\;
(K^*)_{\al}^{\;\; \oga}h_{\oga \be}=K_{\al \be}
\end{eqnarray*}

\noindent and

\begin{eqnarray*}
(K^*)_{\al}^{\;\; \ga}=K_{\al \obe}h^{\ga \obe} \; , \;\;
(K^*)_{\al}^{\;\; \oga}=K_{\al \be}h^{\oga \be}.
\end{eqnarray*}

\noindent Since $K^* (T_{\al}) \in H(M) \otimes \C$, we have
$(K^*)_{\al}^{\;\; 0}=0$, so that

\[
K^*(T_{\al}) = K_{\al \obe}h^{\ga \obe}T_{\ga} + K_{\al \be}h^{\oga
\be}T_{\oga}\, .
\]

\noindent Similarly,

\[
K^*(T_{\oal}) = K_{\oal \obe}h^{\ga \obe}T_{\ga} + K_{\oal \be}h^{\oga
\be}T_{\oga} = \overline {K^*(T_{\al})}\, .
\]

\noindent Sometimes we will identify $K$ with its adjoint
representation $K^*$. With this identification, we say that the
tensor $K$ of type $(2,0)$ can be seen as a tensor of type $(1,1)$ and
denote it by $(K^*)_{A}^{\;\; B} = K_{A}^{\;\; B}$. The matrix
representation of $K^*$ with respect to a frame
$\{T_{\al},T_{\oal},T\}$ is

\[
K^*: \begin{bmatrix}
K_{\al}^{\;\; \ga} & K_{\oal}^{\;\; \ga} & 0\\
K_{\al}^{\;\; \oga} & K_{\oal}^{\;\; \oga} & 0\\
0 & 0 & 0
\end{bmatrix}.
\]

\noindent If $\{T_{\al}\}$ is such that $h_{\al \obe} = \delta_{\al \be}$, we derive

\[
K_{\al}^{\;\; \ga} = K_{\al \obe}h^{\ga \obe} = K_{\al \oga}\, .
\]

\noindent Since $K_{\al \oga}$ is a hermitian matrix, we then deduce that
$K_{\al}^{\,\,\, \ga}$ is a hermitian matrix. In what follows,
we focus our attention on the hermitian block $K_{\al}^{\;\;
\ga}$ of the adjoint representation of $K$. It is now appropriate to
define invariants related to hermitian matrices.

Let $A$ be a hermitian matrix of order $n$. For each $k \in
\{1,...,n\}$, we define the \textbf{$k$-th invariant $\sigma_k(A)$
associated to the matrix $A$} as the $k$-th symmetric elementary
function of the eigenvalues of $A$. Precisely, if $\{\la_1, \ldots, \la_n\}$ are the eigenvalues of $A$, we have

\[
\sigma_k(A) = \sum_{i_1 < \ldots < i_k} \lambda_{i_1} \ldots
\lambda_{i_k}\, .
\]

\n Obviously, $\sigma_k(U^*AU)=\sigma_k(A)$ for any unit matrix $U$, $\sigma_1(A)=\mbox{tra\c{c}o}(A)$
and $\sigma_n(A)=\mbox{det}(A)$.

For each $k \in \{1,...,n\}$, consider the G\"{a}rding's cone in $\R^n$ defined by

\[
\Gamma_k^+ = \{\Lambda = (\lambda_1, \ldots, \lambda_n) \in \R^n :\; \sigma_j(\Lambda) > 0\ \mbox{for all}\ j \leq k\}\ .
\]

\n A hermitian matrix $A$ of order $n$ is said to be $k$-positive, if the $n$-tuple $\Lambda = (\lambda_1, \ldots, \lambda_n)$ of eigenvalues of $A$ belongs to $\Gamma_k^+$. In other words, $\sigma_j(A) > 0$ for all $j \leq k$. In this case, we denote $A \in \Gamma_k^+$.

Let us turn our attention to the pseudohermitian Schouten tensor
$S_{\ta}$. The adjoint representation of $S_{\ta}$, with respect to a
frame $\{ T_{\al} \}$, is given by

\[
S_{\ta}^*: \begin{bmatrix}
S_{\al}^{\;\; \ga} & S_{\oal}^{\;\; \ga} & 0\\
S_{\al}^{\;\; \oga} & S_{\oal}^{\;\; \oga} & 0\\
0 & 0 & 0
\end{bmatrix},
\]

\noindent where

\[
S_{\al}^{\;\; \ga}=S_{\al \obe}h^{\ga \obe}=\frac{1}{n+2} ( {
R_{\al}^{\;\; \ga}-\frac{1}{2(n+1)}R \dl^{\ga}_{\al}} )\, .
\]

\noindent Remark that

\begin{eqnarray*}
\sigma_1(S_{\al}^{\;\; \ga})=S_{\al}^{\;\; \al}=\frac{1}{n+2} ( {
R_{\al}^{\;\; \al}-\frac{1}{2(n+1)}R \dl^{\al}_{\al}})=
\frac{1}{n+2} { \frac{2(n+1)R-Rn}{2(n+1)}}\, ,
\end{eqnarray*}

\noindent so that

\[
\sigma_1(S_{\al}^{\;\; \ga})=\frac{R}{2(n+1)}\, .
\]

\noindent So, the $k$-th invariant $\sigma_k$ of $S_{\al}^{\;\; \ga}$ is a good generalization of the Webster scalar
curvature $R$ associated to $\ta$. Thus, given a CR manifold $(M, \theta)$ of dimension $2n + 1$, we
define its \textbf{pseudohermitian $k$-curvature} by

\[
\sigma_k(\ta)=\sigma_k(S_{\al}^{\;\; \ga})\, .
\]

A quite useful notion, closely related to the $k$-positive metric one, is the notion of $k$-positive pseudohermitian structure. In fact, a pseudohermitian structure $\ta$ on $M$ is said to be $k$-positive, if $\sigma_j(\ta)(x) > 0$ for all $x \in M$ and $1 \leq j \leq k$. When both $\tilde{\ta} = e^{2u} \ta$ and $\ta$ are $k$-positive, we simply say that $u$ is $k$-admissible.

Another important geometric tensor is the pseudohermitian Cotton tensor, which can be defined from the Schouten
tensor. Precisely, given a CR manifold $(M, \ta)$ of dimension $2n +
1$, the \textbf{pseudohermitian Cotton tensor} associated to $\theta$ is the tensor $C_\ta$ of type $(3,0)$ locally defined by

\begin{eqnarray}
C_{\rro \osi \ga}=S_{\rro \osi ; \ga}-S_{\ga \osi ; \rro}\, ,
\end{eqnarray}

\noindent where again $S_{\rro \osi}$ denote the components of the pseudohermitian Schouten tensor $S_{\ta}$. Note that this tensor measures a certain
symmetry of the covariant derivatives of the pseudohermitian Schouten tensor such as occurs for the Riemannian Cotton tensor.

In \cite{Via1}, Viaclovsky proved that the $k$-Yamabe problem is variational whenever the Cotton tensor of the metric $g$ is null. Theorem 1.1 establishes the corresponding result in the CR geometry context. A sufficient condition for the pseudohermitian Cotton tensor to be null is the pseudohermitian Schouten tensor to be parallel with respect to the Tanaka-Webster connection, that is $\nabla S_\theta = 0$. Another sufficient condition, which includes in particular the pseudohermitian spatial forms presented by Webster in \cite{Web}, is the pseudohermitian torsion to be parallel with respect to this same connection, this is $\nabla \tau_\theta = 0$.

In order to introduce the Cotton-admissibility notion, we remark that the condition $C_{\tta}=0$ is not CR invariant. In fact, denoting by $\til{C}$ and $\til{S}$ and $C$ and $S$, respectively, the pseudohermitian Cotton and Schouten tensors associated to $\tta = e^{2u} \ta$ and $\ta$, we find

\begin{eqnarray} \label{tilcotten} \til C_{\al \obe ; \si}=C_{\al \obe ;
\si}-2iu^{\rro}(h_{\si \obe}A_{\al \rro}+h_{\al \obe}A_{\si
\rro})-2(u_{\al}u_{\obe \si}-u_{\si}u_{\obe \al}).
\end{eqnarray}

\n Indeed,

\begin{eqnarray*}
\til C_{\al \obe ; \si}&=&\til S_{\al \obe ; \si}-\til S_{\si \obe ; \al}\\
&=&(S_{\al \obe} -u_{\al \obe} - u_{\obe \al} -2u_{\ga}u^{\ga}h_{\al \obe} )_{; \si}-(S_{\si \obe} -u_{\si \obe} - u_{\obe \si} -2u_{\ga}u^{\ga}h_{\si \obe} )_{ ; \al}\\
&=& S_{\al \obe  ; \si} -u_{\al \obe  ; \si} - u_{\obe \al  ; \si} -2u_{\ga  ; \si}u^{\ga}h_{\al \obe} -2u_{\ga}u_{\;\; ; \si}^{\ga}h_{\al \obe}\\
&&-(S_{\si \obe ; \al} -u_{\si \obe; \al} - u_{\obe \si; \al}
-2u_{\ga; \al}u^{\ga}h_{\si \obe}-2u_{\ga}u_{\;\; ; \al}^{\ga}h_{\si
\obe} )\, ,
\end{eqnarray*}

\n so that

\begin{eqnarray*}
\til C_{\al \obe ; \si}&=&C_{\al \obe ; \si}-(u_{\al \obe  ; \si}-u_{\si \obe; \al}) - (u_{\obe \al  ; \si}- u_{\obe \si; \al})\\
&&-2(u_{\obe}u_{\ga \si}+u_{\ga}u_{\obe \si})h^{\ga \obe}h_{\al \obe}+2(u_{\obe}u_{\ga \al}+u_{\ga}u_{\obe \al})h^{\ga \obe}h_{\si \obe}\\
&=&C_{\al \obe ; \si}-(u_{\al \obe  ; \si}-u_{\si \obe; \al})-(u_{\obe \al  ; \si}- u_{\obe \si; \al})\\
&&-2(u_{\obe}u_{\ga \si}+u_{\ga}u_{\obe \si})\dl^{\;\; \ga}_{\al}+2(u_{\obe}u_{\ga \al}+u_{\ga}u_{\obe \al})\dl^{\;\; \ga}_{\si}\\
&=&C_{\al \obe ; \si}-(u_{\al \obe  ; \si}-u_{\si \obe; \al})-(u_{\obe \al  ; \si}- u_{\obe \si; \al})\\
&&-2(u_{\obe}u_{\al \si}+u_{\al}u_{\obe \si})+2(u_{\obe}u_{\si \al}+u_{\si}u_{\obe \al})\\
&=&C_{\al \obe ; \si}- (u_{\al \obe  ; \si}-u_{\si \obe; \al}) - (u_{\obe \al  ; \si}- u_{\obe \si; \al}) + 2(u_{\al}u_{\obe
\si}-u_{\si}u_{\obe \al}) \, .
\end{eqnarray*}

\n Thus, the desired equality follows directly from the identities (see \cite{BD1}, \cite{Lee2})

\[
u_{\al \obe  ; \si}-u_{\si \obe; \al} = iu^{\rro}(h_{\si \obe}A_{\al \rro}-h_{\al \obe}A_{\si \rro})\, ,
\]

\[
u_{\obe \al  ; \si}- u_{\obe \si; \al} = iu^{\rro}(h_{\al \obe}A_{\si \rro} - h_{\si \obe}A_{\al \rro})\, .
\]

\n Now, if the pseudohermitian torsion of $\ta$ is parallel with respect to the Tanaka-Webster connection, then $C_{\al \obe \ga}=0$. So, from (\ref{tilcotten}), we deduce that

\[
\til C_{\al \obe ; \si}=-2iu^{\rro}(h_{\si \obe}A_{\al \rro} + h_{\al \obe}A_{\si \rro})-2(u_{\al}u_{\obe \si}-u_{\si}u_{\obe \al})\, .
\]

\n On the other hand, if the pseudohermitian torsion is null, we obtain

\[
\til C_{\al \obe ; \si}=-2(u_{\al}u_{\obe \si}-u_{\si}u_{\obe \al}).
\]

\n This suggests the following definition. A function $u \in C^{\infty}(M)$ is said to be \textbf{Cotton-admissible}, if it satisfies

\[
u_{\al}u_{\obe \si} - u_{\si}u_{\obe \al} = 0
\]
for any $\al , \be , \si \in \{1, ... ,n\}$. Some examples of Cotton-admissible functions $u \in C^{\infty}(M)$ are given by CR-pluriharmonic functions, that is functions $u \in C^\infty(M)$ for which there is a function $v \in C^\infty(M)$ such that $f=u+iv$ is CR-holomorphic. Thus, if $u$ is CR-pluriharmonic and $\{ T_{\al}\}$ is a frame of $(M, \ta)$, then

\[
T_{\obe}f = T_{\obe}u+iT_{\obe}v = u_{\obe}+iv_{\obe}=0\, .
\]

\n Hence, $u_{\obe} = 0$ and, therefore, $u_{\obe \al} = u_{\obe \si}=0$ for any $\al ,\si \in \{1, ... ,n\}$.

In the next sections, we need the following sets

\[
[\theta]_{+} = \{ e^{2u}\ta : u \mbox{ is $k$-admissible} \}\, ,
\]

\[
\mathcal{C}[\ta] = \{ e^{2u}\ta : u \mbox{ is Cotton-admissible} \}
\]

\n and

\[
\mathcal{C}[\ta]_{+} = \mathcal{C}[\ta] \cap [\theta]_{+} = \{ e^{2u}\ta:u \; \mbox{is $k$-admissible and Cotton-admissible} \}\, .
\]

\n A pseudohermitian structure $\tta \in \mathcal{C}[\ta]$ will be called Cotton-admissible.\\

\section{The CR $k$-Yamabe problem}

Once introduced the notion of pseudohermitian $k$-curvature in the previous section, it naturally arises the following question:

\begin{center}
\textbf{The CR $k$-Yamabe Problem.} Given a compact CR manifold $(M, \theta)$ of dimension $2n+1 \geq 3$, finding an
orientation preserving pseudohermitian structure conformal to $\theta$ with constant pseudohermitian $k$-curvature.
\end{center}

Let $(M, \ta)$ be a compact CR manifold of dimension $2n + 1$. Given a pseudohermitian structure $\tta = e^{2u} \ta$, according to the relation (\ref{schout}), the Schouten tensors $S_{\tta}$ and $S_{\ta}$ associated to $\tta$ and $\ta$ relate in an admissible coframe as

\[
S_{\tta} = S_{\ta} - 2 \grad^{2} u + (i T_\ta u-\normc {du} 2 \ta ) L_{\ta}\, .
\]

Assume that $\tta$ has constant pseudohermitian $k$-curvature, that is $\si_k(\tta) = \la$ for some constant $\la$. Since

\begin{eqnarray*}
\si_k(\tta)&=&\si_k({\til S}_{\al}^{\;\; \ga})=\si_k(\tih^{\ga \obe}{\til S}_{\al \obe})\\
&=&e^{-2ku}\si_k(h^{\ga \obe}{\til S}_{\al \obe})\\
&=&e^{-2ku}\si_k(S_{\al}^{\;\; \ga} -2u_{\al}^{\;\; \ga}+(i T_\ta u-\normc {du} 2 \ta )\dl_{\al}^{\;\; \ga})\, ,
\end{eqnarray*}

\n it follows that $\si_k(\tta) = \la$ if, and only if, $u$ satisfies the equation

\[
\si_k(S_\ta -2\grad^{2}u+(i T_\ta u-\normc {du} 2 \ta )L_{\ta}) = \la e^{2ku}\, .
\]

We next show that the preceding equation yields the CR Yamabe equation when $k = 1$. In particular, this equation will be called the \textbf{CR $k$-Yamabe equation}. In fact, we have

\begin{eqnarray*}
\si_1\lp S_{\ta} -2\grad^{2}u+(i T_\ta u - \normc {du} 2 \ta )L_{\ta} \rp&=&\mbox{trace} (  (S_{\al \obe} -2u_{\al \obe}+(iu_0 - \normc {du}
2 \ta) h_{\al \obe}) h^{\ga \obe} )\\
&=&\mbox{trace} (S_{\al \obe} -2u_{\al \obe}+iu_0 h^{\ga \obe} - \normc {du} 2 \ta h_{\al \obe}h^{\ga \obe} )\\
&=&\mbox{trace} (  S_{\al}^{\;\; \ga} -u_{\al}^{\;\; \ga}-u^{\ga}_{\;\; \al}- \normc {du} 2 \ta \dl_{\al}^{\;\; \ga} )\\
&=&  S_{\al}^{\;\; \al} -u_{\al}^{\;\; \al}-u_{\oal}^{\;\;\; \oal}- \normc {du} 2 \ta \dl_{\al}^{\;\; \al} \\
&=& \frac{R}{2(n+1)} +\lap_b u - n\normc {du} 2 \ta \, ,
\end{eqnarray*}

\n where $\lap_b$ denotes the well-known sublaplacian operator. Here, it is used that $\lap_b$ possesses a particularly simple expression
in terms of covariant derivatives (see \cite{Lee1}) when one considers an admissible coframe $\{ \ta^{\al} \}$, namely

\begin{equation} \label{sublaplaciano}
\lap_b u = -(u_{\al}^{\;\; \al} + u_{\oal}^{\;\; \oal})\, .
\end{equation}

\n Thus,

\[
R + 2(n+1)\lap_b u - 2n(n+1)\normc {du} 2 \ta = 2(n+1) \la e^{2u}\, .
\]

\n Performing the canonical change $e^{2u}=v^{p-2}$, with $p = 2 + 2/n$, modulo a constant factor, it arrives at the CR Yamabe equation

\begin{equation} \label{CR}
{\cal L}_1[v]:= p \, \lap_b  v + R v = 2(n+1) \la v^{p-1}\, .
\end{equation}

\n The operator ${\cal L}_1$ is the well-known CR Yamabe conformal operator.

As mentioned in the introduction, the CR $k$-Yamabe problem was
completely solved in the case that $k = 1$ by Jerison, Lee, Gamarra
and Yacoub based on the study of the above equation. Indeed, they
approached and solved the problem under a variational view, since
the equation (\ref{CR}) can be seen as the Euler-Lagrange equation
associated to the minimization problem

\[
\la(M) = \inf\{Y(v):\; G(v)=1 \}\, ,
\]

\n where

\[
Y(v):=\int_M v {\cal L}_1[v] \; dV_{\ta}=\int_M (p \, \normc {dv} 2 \ta + Rv^2) \; dV_{\ta}
\]

\n and

\[
G(v):= \int_M dV_{\tta}=\int_M |v|^p\; dV_{\ta}\, .
\]

\n One of the main results of this work concerns with the variational nature of the CR $k$-Yamabe problem when $k > 1$. In particular, we apply it in order to investigate the problem on the CR sphere.

Dealing specifically with the proposed problem, our main contributions can be resumed in the two results below, which will be proved in Sections 4 and 5. Consider the following generalization of the functional $Y$:

\[
Y_{k}: \tilde \ta \mapsto \int_{M} \si_k(\tilde \ta) \, dV_{\tilde \ta}
\]

\n constrained to the subset

\[
{\mathcal M}=\{\tta \in [\ta]: \int_M dV_{\tta} =1\}
\]

\n of conformal pseudohermitian structures $\tta$ of unit volume and compatible with $\theta$.

\begin{teorema} (The variational property)
Let $(M, \ta)$ be a compact CR manifold of dimension $2n+1$. If $k \geq 1$, then a pseudohermitian structure $\tilde{\ta} \in {\cal M}$ conformal to $\ta$ such that $C_{\tilde \ta} = 0$ is a critical point of $Y_k$ if, and only if, $\sigma_k (\tilde \ta) = \lambda_k$ for some constant $\lambda_k$.
\end{teorema}

\begin{teorema} (An Obata type classification result)
Let $(\s^{2n+1}, \hat{\ta})$ be the CR sphere of dimension $2n + 1$ and be $k \geq 1$. If the pseudohermitian structure $\tilde{\ta} = e^{2u} \hat{\ta} \in \mathcal{C}[\hat{\ta}]_{+}$ has constant pseudohermitian $k$-curvature, then $\tilde{\ta}$ is equal to a multiple of the canonical form $\hat{\ta}$ by some CR-automorphism of $\s^{2n+1}$.
\end{teorema}

\section{The variational property and pseudo-Einstein manifolds}

Our main goal in this section is to provide the proof of Theorem 1.1 (or Theorem 3.1). For this, we need one preliminary lemma. Before stating it, we recall that the \textbf{$k$-th Newton transformation} associated to a hermitian matrix $A$ is given by

\[
T_k(A)=\sigma_k(A)I-\sigma_{k-1}(A)A + \ldots + (-1)^k A^k\, .
\]

\n Some basic properties are satisfied by this transformation. For example, using the fact that if $A$ and $B$ are hermitian matrices, then $AB$ is a
hermitian matrix if, and only if, $AB=BA$, one concludes that $T_k(A)$ is always a hermitian matrix.

Let $1 \leq k \leq n$, $1 \leq i_1, \ldots, i_k \leq n$, the generalized Kronecker
symbol is defined by

\[
\dl_{j_1 \ldots j_k}^{i_1 \ldots i_k}=\left\{\begin{array}{rl}
1,  & \mbox{if } i_1, \ldots , i_k \mbox{ are distinct and } (j_1, \ldots, j_k) \mbox{ is an even permutation of }(i_1, \ldots, i_k), \\
-1, & \mbox{if } i_1, \ldots, i_k \mbox{ are distinct and } (j_1, \ldots, j_k) \mbox{ is an odd permutation of }(i_1, \ldots, i_k), \\
0,  & \mbox{otherwise}
\end{array}\right.
\]

In the following lemma (see \cite{Rei}), we summarize some of the main properties satisfied by $\sigma_k(A)$ and $T_k(A)$:

\begin{lema}
\label{Sigmak} Let $A=\lp A_{i}^{\; j}\rp$ be a hermitian matrix of order $n$.
Then:

\[
\sigma_k(A)=\frac{1}{k!}\sum{\dl_{j_1 \ldots j_k}^{i_1 \ldots i_k}}A_{i_1}^{\;\; j_1} \ldots A_{i_k}^{\;\; j_k}\, ,
\]

\[
{T_k(A)}_{j}^{\;\;\; i}=\frac{1}{k!}\sum{\dl_{j_1 \ldots j_k j}^{i_1 \ldots i_k i}}A_{i_1}^{\;\; j_1} \ldots A_{i_k}^{\;\; j_k}\, ,
\]

\[
\sigma_1(T_k(A) \circ A)=(k+1)\sigma_{k+1}(A)\, ,
\]

\[
T_k(A)=\sigma_k (A)I-T_{k-1}(A)A\, ,
\]

\[
\sigma_1(T_k(A))=(n-k)\sigma_k(A)\, .
\]

\end{lema}

We now are ready to prove Theorem 1.1.

\begin{proof}
We first write $\tta= e^{2u} \ta$ and evoke the relation (\ref{schout}), between the Schouten tensors $S_{\tta}$ and $S_\ta$, so that

\[
\si_k(\tta) = e^{-2ku}\si_k(S_\ta -2\grad^{2}u+(i T_\ta u-\normc {du} 2 \ta )L_{\ta})\, .
\]

\n Consider now the functional $F_k$ on $C^{\infty}(M)$ given by

\[
F_k(u) : = Y_k(e^{2u} \ta) = \int_M e^{-2(n+k+1)u}\si_k(S_\ta -2\grad^{2}u+(i T_\ta u-\normc {du} 2
\ta )L_{\ta})\; dV_{\ta}\, .
\]

\n Here it is used that $dV_{\tta} = e^{-2(n+1)u}dV_{\ta}$.

Let $u, \phi \in C^{\infty}(M)$ and $u(t)$ be a curve in $C^{\infty}(M)$ such that $u(0)=u$ and $\dot{u}(0)=\phi$. A simple computation gives us

\[
\left.\frac{d}{dt}\right|_{t=0} F_k(u(t))=-2(n+k+1)\int_M \phi e^{-2(n+k+1)u}\si_k(S_\ta -2\grad^{2}u+(i T_\ta u-\normc {du} 2 \ta )L_{\ta})\; dV_{\ta}
\]

\[
+\int_M e^{-2(n+k+1)u} \left.\frac{d}{dt}\right|_{t=0} \si_k(S_{\tta})\; dV_{\ta}\, .
\]

\n Using the two first relations of Lemma \ref{Sigmak}, by direct differentiation, we easily deduce that

\[
\frac{d}{dt} \si_k(A(t))=T_{k-1}(A(t))_{\be}^{\;\; \al}\frac{d}{dt}
A(t)_{\al}^{\;\; \be}=\si_1 ( {T_{k-1}(A(t))\circ \frac{d}{dt}
A(t)} )\, .
\]

\n So, we obtain

\begin{eqnarray*}
\left.\frac{d}{dt}\right|_{t=0} \si_k(S_{\tta}) &=& \left.\frac{d}{dt}\right|_{t=0} \si_k \lp S_\ta -2\grad^{2}u(t)+\lp i T_\ta u(t)-L_{\ta}^{*}(du(t),du(t)) \rp L_{\ta}\rp\\
&=& \si_1 ({T_{k-1}(S_{\tta})\circ \frac{d}{dt} (S -2\grad^{2}u(t)+(i T_\ta u(t)-L_{\ta}^{*}(du(t),du(t)) )L_{\ta})})\\
&=& \si_1 ( {e^{(k-1)u} T_{k-1}(S_{\tta}) \circ (-2\grad^{2}\phi+(i T_\ta \phi-L_{\ta}^{*}(du,d\phi)-L_{\ta}^{*}(d\phi,du) )L_{\ta})})\, .
\end{eqnarray*}

\n On the other hand, we have

\[
(-2\grad^{2}\phi+(i T_\ta \phi-L_{\ta}^{*}(du,d\phi)-L_{\ta}^{*}(d\phi,du)
)L_{\ta})_{\al \obe} =(-2\phi_{\al \obe}+i\phi_0 h_{\al \obe})-(L_{\ta}^{*}(du,d\phi)-L_{\ta}^{*}(d\phi,du) )h_{\al \obe}\, .
\]

\n Thus, since

\[
-2\phi_{\al \obe}+i\phi_0 h_{\al \obe} = -\phi_{\al \obe}-\phi_{\obe \al}
\]

\n and

\[
L_{\ta}^{*}(du,d\phi) = u_{\al}\phi_{\obe}h^{\al \obe}+u_{\oal}\phi_{\be}h^{\be \oal} = u^{\obe}\phi_{\obe}+u^{\be}\phi_{\be}=L_{\ta}^{*}(d\phi,du)\, ,
\]

\n we derive

\[
(-2\grad^{2} \phi + (i T_\ta \phi-L_{\ta}^{*}(du,d\phi) - L_{\ta}^{*}(d\phi,du) )L_{\ta})_{\al \obe} = -\phi_{\al \obe} - \phi_{\obe \al}-(2u^{\ga}\phi_{\ga} + 2u^{\oga}\phi_{\oga})h_{\al \obe}
\]

\[
= -(\phi_{\al \obe} + 2u^{\ga}\phi_{\ga}h_{\al \obe}) - (\phi_{\obe \al} + 2u^{\oga}\phi_{\oga}h_{\al \obe})\, .
\]

\n From the relation between the Christoffel symbols associated to $\tta$ and $\ta$,

\[
{\til \Ga}^{\ga}_{\obe \al} = \Ga^{\ga}_{\obe \al}-2u^{\ga}h_{\al \obe}\, ,
\]

\n we obtain

\[
(\til \grad^2 \phi)_{\al \obe} = T_{\obe}T_{\al}\phi-{\til \Ga}^{\ga}_{\obe \al}T_{\ga}\phi=(\grad^2 \phi)_{\al
\obe}+2u^{\ga}h_{\al \obe}=\phi_{\al \obe}+2u^{\ga}\phi_{\ga}h_{\al \obe}\, .
\]

\n Similarly,

\[
(\til \grad^2 \phi)_{\obe \al} = \phi_{\obe \al}+2u^{\oga}\phi_{\oga}h_{\al \obe}\, .
\]

\n Therefore,

\[
\left.\frac{d}{dt}\right|_{t=0} \si_k(S_{\tta}) = \si_1 ( {e^{2(k-1)u} T_{k-1}(S_{\tta}) \circ (-e^{2u}((\til \grad^2 \phi)_{\al}^{\;\; \ga} + (\til \grad^2 \phi)^{\ga}_{\;\; \al} ))} )
\]

\[
=-e^{2ku}\si_1 ( { T_{k-1}(S_{\tta}) \circ ((\til \grad^2 \phi)_{\al}^{\;\; \ga} + (\til \grad^2 \phi)^{\ga}_{\;\; \al} )} )\, .
\]

\n Finally, from

\begin{eqnarray*}
(\til \grad^2 \phi)_{\al}^{\;\; \ga} + (\til \grad^2 \phi)^{\ga}_{\;\; \al} &=&((\til \grad^2 \phi)_{\al \obe}+ (\til \grad^2 \phi)_{\obe \al} )\tih^{\ga \obe}\\
&=&(2(\til \grad^2 \phi)_{\al \obe}-i T_{\tta} \phi \tih_{\al \obe})\tih^{\ga \obe}\\
&=&2(\til \grad^2 \phi)_{\al}^{\;\; \ga}-i T_{\tta} \phi \,
\dl_{\al}^{\;\; \ga}\, ,
\end{eqnarray*}

\n we find

\[
\left.\frac{d}{dt}\right|_{t=0} \si_k(S_{\tta}) = -e^{2ku}\si_1 ( { T_{k-1}(S_{\tta}) \circ ( 2\til \grad^2 \phi-i T_{\tta} \phi I  )}) \, .
\]

\n Consequently,

\begin{eqnarray*}
\left.\frac{d}{dt}\right|_{t=0} F_k(u(t)) &=& -2(n+k+1)\int_M \phi e^{-2(n+k+1)u}\si_k(S -2\grad^{2}u+(i T_\ta u-\normc {du} 2 \ta )L_{\ta})\; dV_{\ta}\\
&& +\int_M e^{-2(n+k+1)u} \left.\frac{d}{dt}\right|_{t=0} \si_k(S_{\tta})\; dV_{\ta}\\
&=& -2(n+k+1)\int_M \phi \si_k(S_{\tta})\; dV_{\tta} - \int_M e^{-2(n+k+1)u}e^{2ku}\si_1 ( { T_{k-1}(S_{\tta}) \circ ( 2\til \grad^2 \phi-i T_{\tta} \phi I  )} ) \; dV_{\ta}\\
&=& -2(n+k+1)\int_M \phi \si_k(S_{\tta})\; dV_{\tta} - \int_M \si_1 ( { T_{k-1}(S_{\tta}) \circ ( 2\til \grad^2 \phi-i T_{\tta} \phi I  )} ) \;
dV_{\tta}\, .
\end{eqnarray*}

\n But, since $T_{k-1}(S_{\tta})$ is hermitian, it follows by the integration by parts that

\[
\int_M \si_1 ( { T_{k-1}(S_{\tta}) \circ ( - 2\til \grad^2 \phi + i T_{\tta} \phi I  )} ) \; dV_{\tta} = -2 Re{\int_M {T_{k-1}(S_{\tta})}_{\ga \;\; ; \al}^{\;\; \al}}u^{\ga}\; dV_{\tta}\, ,
\]

\n so that

\[
\left.\frac{d}{dt}\right|_{t=0} F_k(u(t))=-2(n+k+1)\int_M \phi \si_k(S_{\tta})\; dV_{\tta}-2 Re{\int_M {T_{k-1}(S_{\tta})}_{\ga \;\; ;
\al}^{\;\; \al}}u^{\ga}\; dV_{\tta}\, .
\]

\n Using now the assumption that $C_{\tta}=0$, we find

\[
\til S_{\be \;\; ; \al}^{\;\; \ga}-\til S_{\al \;\; ; \be}^{\;\; \ga}=(\til S_{\be \osi ; \al}-\til S_{\al \osi ; \be})\tih^{\ga
\osi}=\til C _{\be \osi \al}\tih^{\ga \osi}=0\, .
\]

\n By Lemma 4.1, we have

\begin{eqnarray*}
{T_{k-1}(S_{\tta})}_{\ga \;\; ; \al}^{\;\; \al}&=&\frac{1}{(k-1)!}\, \dl_{\ga_1 ... \ga_{k-1} \ga}^{\al_1 ... \al_{k-1} \al}(\til S_{\al_1}^{\;\;\; \ga_1} \ldots \til S_{\al_{k-1}}^{\;\;\;\;\;\;\; \ga_{k-1}})_{; \, \al}\\
&=&\frac{1}{(k-2)!} \, \dl_{\ga_1 ... \ga_{k-1} \ga}^{\al_1 ...
\al_{k-1} \al}\til S_{\al_1 \;\; ; \, \al}^{\;\;\; \ga_1} \ldots
\til S_{\al_{k-1}}^{\;\;\;\;\;\;\; \ga_{k-1}}\, .
\end{eqnarray*}

\n So, we are led to

\[
{T_{k-1}(S_{\tta})}_{\ga \;\; ; \al}^{\;\; \al} = 0
\]

\n and this implies that

\[
\left.\frac{d}{dt}\right|_{t=0} F_k(u(t))=-2(n+k+1)\int_M \phi \si_k(S_{\tta})\; dV_{\tta}\, .
\]

\n Since we are considering the functional $Y_k$ constrained to the conformal pseudohermitian structures of unit volume, the remaining conclusion follows readily from Lagrange multipliers.

\end{proof}

Let $(M, \ta)$ be a CR manifold of dimension $2n + 1$. We say that $\ta$ is pseudo-Einstein, if $R_{\al \obe} = \frac{R}{n}h_{\al \obe}$. Unlike the Riemannian case, the pseudo-Einstein condition in general does not imply that the Webster scalar curvature is constant (for instance, see \cite{Lee2}). However, $\si_k(\ta)$ is constant if, and only if, the Webster scalar curvature is constant, whenever $\ta$ is pseudo-Einstein. In fact,

\[
S_{\al \obe}=\frac{1}{n+2}( {R_{\al \obe}-\frac{R}{2(n+1)}}h_{\al
\obe} )=\frac{1}{n+2}( {\frac{1}{n}-\frac{1}{2(n+1)} } )
Rh_{\al \obe}=\frac{R}{2n(n+1)}h_{\al \obe}\, ,
\]

\n so that

\[
\si_k(\ta)=\si_k(S_{\al}^{\;\; \ga})=\si_k ( \frac{R}{2n(n+1)}
\dl_{\al}^{\;\; \ga} )=C(n,k)\frac{R^k}{{\lp 2n(n+1) \rp}^k}\, ,
\]

\n where $C(n,k)$ denotes the Newton binomial coefficient.

A natural question that arises here is whether the solutions of the CR $k$-Yamabe problem
on pseudo-Einstein CR manifolds coincide with the
solutions of the CR Yamabe problem. In an attempt to answer this question, we are led to analyze the effect of the deformation $\tta =
e^{-2u}\ta$ of a pseudo-Einstein structure $\ta$. In
\cite{Lee2}, Lee showed that if $\ta$ is pseudo-Einstein,
a necessary and sufficient condition for $\tta$ to be pseudo-Einstein
is $u$ to be CR-pluriharmonic. So, if $\ta$ is
pseudo-Einstein and if $\tta=e^{-2u}\ta$ is a conformal deformation of
$\ta$ by a CR-pluriharmonic function $u$, then $\tta$ is a solution of the CR Yamabe problem if, and only if, it is a solution of the CR $k$-Yamabe problem.

As an immediate consequence of Theorem 2.1, it follows that if $\tta$ is a solution of the $k$-Yamabe problem and, in addition, the Cotton tensor
$C_{\tta}$ vanishes, then $\tta$ is pseudo-Einstein. In particular, $\si_1(\tta)$ is constant and, consequently, solutions of the CR $k$-Yamabe problem are also solutions of the CR Yamabe problem.

In order to prove this fact, beyond Lemma 4.1, we need another key tool, see \cite{Via1}.

\begin{lema}
\label{Lema2} Let $A$ be a hermitian matrix of order $n$. If $A \in
\Ga^{+}_{k}$, then

\[
\si_{k+1}(A) \leq \frac{n-k}{n(k+1)}\si_k(A)\si_1(A)\, .
\]

\n Furthermore, $$\si_{k+1}(A) = \frac{n-k}{n(k+1)}\si_k(A)\si_1(A)$$
if, and only if, $A=\la I_n$ for some constant $\la$, where $I_n$ stands for the identity matrix of order $n$.
\end{lema}

We now are ready to provide a proof of Theorem 1.2.

\begin{proof}
We begin by writing $\tta = e^{-2u}\ta$. Since $\ta$ is pseudo-Einstein, we have

\[
S_{\ta} = \frac{R}{2n(n+1)}L_{\ta}\, .
\]

\n Using the identities

\[
R=e^{-2u}({\til R}+2(n+1)(\til{\lap}_b u)-2n(n+1)\|du\|^{2}_{\tta})
\]

\n and

\[
L_{\ta} = e^{2u} L_{\tta}\, ,
\]

\n we deduce that

\[
S_\ta = ( \frac{\til{R}}{2n(n+1)}+\frac{\til{\lap}_b u}{n}-\|du\|^{2}_{\tta} ) L_{\tta}\, .
\]

\n Then, joining the above equality with

\[
S_\ta = S_{\tta} -2\til{\grad}^2 u +(iT_{\tta}u - \|du\|^{2}_{\tta})L_{\tta}\, ,
\]

\n we obtain

\[
2\til{\grad}^2 u = S_{\tta} - \dfrac{\si_1(S_{\tta})}{n}L_{\tta}+ ( iT_{\tta}u - \dfrac{\til{\lap}_b
u}{n} ) L_{\tta}\, .
\]

\n So, we can write

\begin{eqnarray*}
\int_M \si_1(T_k(S_{\tta}) \circ (2\til{\grad}^2 u))\; dV_{\tta} &=& \int_M \si_1 ( T_k(S_{\tta}) \circ ( S_{\tta} - \dfrac{\si_1(S_{\tta})}{n} I + ( iT_{\tta}u - \dfrac{\til{\lap}_b u}{n} ) I ) )\; dV_{\tta} \\
&=&\int_M \si_1 ( T_k(S_{\tta}) \circ S_{\tta} - \dfrac{\si_1(S_{\tta})}{n} T_k(S_{\tta}) + ( iT_{\tta}u - \dfrac{\til{\lap}_b u}{n} ) T_k(S_{\tta}) )\; dV_{\tta}\\
&=&\int_M \si_1 ( T_k(S_{\tta}) \circ S_{\tta}) - \dfrac{\si_1(S_{\tta})}{n} \si_1( T_k(S_{\tta}) ) + ( iT_{\tta}u - \dfrac{\til{\lap}_b u}{n} ) \si_1 ( T_k(S_{\tta}) )\; dV_{\tta}\, .
\end{eqnarray*}

\n By Lemma \ref{Sigmak}), we have

\[
\si_1(T_k(S_{\tta}) \circ S_{\tta})=(k+1)\si_{k+1}(S_{\tta}) \;\;\; \mbox{and} \;\;\; \si_1(T_k(S_{\tta}))=(n-k)\si_k(S_{\tta})\, ,
\]

\n so that

\[
\int_M \si_1(T_k(S_{\tta}) \circ (2\til{\grad}^2 u))\; dV_{\tta} =  \int_M (k+1)\si_{k+1}(S_{\tta}) - \dfrac{\si_1(S_{\tta})}{n} (n-k)\si_k(S_{\tta}) +  ( iT_{\tta}u - \dfrac{\til{\lap}_b u}{n} )   (n-k)\si_k(S_{\tta}) \; dV_{\tta}
\]

\[
= (k+1) \int_M \si_{k+1}(S_{\tta}) - \dfrac{n-k}{n(k+1)}\si_k(S_{\tta}) \si_1(S_{\tta}) \; dV_{\tta} + \int_M iT_{\tta} u(n-k)\si_k(S_{\tta}) \; dV_{\tta} - \dfrac{1}{n}\int_M \til{\lap}_b u \si_k(S_{\tta})\; dV_{\tta}\, .
\]

\n Consequently,

\[
\int_M \si_1(T_k(S_{\tta}) \circ (2\til{\grad}^2 u))\; dV_{\tta} - \int_M iT_{\tta} u \si_1(T_k(S_{\tta}))\; dV_{\tta}
\]

\[
= (k+1) \int_M \si_{k+1}(S_{\tta}) - \dfrac{n-k}{n(k+1)}\si_k(S_{\tta}) \si_1(S_{\tta}) \; dV_{\tta} -\dfrac{1}{n}\int_M \til{\lap}_b u \si_k(S_{\tta})\; dV_{\tta}\, .
\]

\n Since  $\si_k(S_{\tta})$ is constant, we have

\[
\int_M \til{\lap}_b u \si_k(S_{\tta})\; dV_{\tta} = \si_k(S_{\tta}) \int_M \til{\lap}_b u\; dV_{\tta} = 0\, .
\]

\n Noting also that

\[
\int_M \si_1(T_k(S_{\tta}) \circ (2\til{\grad}^2 u))\; dV_{\tta} - \int_M iT_{\tta} u \si_1(T_k(S_{\tta}))\; dV_{\tta} =  \int_M \si_1(T_k(S_{\tta}) \circ (2 \til{\grad}^2 u -iT_{\tta} u I))\; dV_{\tta}
\]

\n and integrating by parts, we find

\[
\int_M \si_1(T_k(S_{\tta}) \circ (2 \til{\grad}^2 u -iT_{\tta} u I))\; dV_{\tta} = 2 Re{\int_M {T_{k-1}(S_{\tta})}_{\ga \;\; ; \al}^{\;\;
\al}}u^{\ga}\; dV_{\tta}\, .
\]

\n But, as previously remarked, the assumption $C_{\tta}=0$ yields ${T_{k-1}(S_{\tta})}_{\ga \;\; ; \al}^{\;\; \al} = 0$. Therefore,

\[
\int_M \si_{k+1}(S_{\tta}) - \dfrac{n-k}{n(k+1)}\si_k(S_{\tta}) \si_1(S_{\tta})\; dV_{\tta} = 0\, .
\]

\n Since $\tta$ is $k$-positive, we can apply the Lemma \ref{Lema2}, so that

\[
\si_{k+1}(S_{\tta}) - \dfrac{n-k}{n(k+1)}\si_k(S_{\tta}) \si_1(S_{\tta}) \leq 0\, .
\]

\n Hence,

\[
\si_{k+1}(S_{\tta}) - \dfrac{n-k}{n(k+1)}\si_k(S_{\tta}) \si_1(S_{\tta}) = 0
\]

\n and then $S_{\tta} = \la I_n$ for some constant $\la$. So,

\[
\si_1(S_{\tta}) = \dfrac{\til{R}}{2(n+1)} = n \la
\]

\n and

\[
\la = \dfrac{\til{R}}{2n(n+1)}\, .
\]

\n Finally, we show that $\tta$ is pseudo-Einstein. Precisely,

\begin{eqnarray*}
\til{R}_{\al \obe}&=&(n+2) \til{S}_{\al \obe}+\frac{\til{R}}{2(n+1)} \til{h}_{\al \obe}\\
&=&(n+2)(\til{S}_{\al}^{\;\; \ga}) \til{h}_{\ga \obe} + \frac{\til{R}}{2(n+1)} \til{h}_{\al \obe}\\
&=&(n+2)( \dfrac{\til{R}}{2n(n+1)}\dl_{\al}^{\;\; \ga} ) \til{h}_{\ga \obe} + \dfrac{\til{R}}{2(n+1)} \til{h}_{\al \obe}\\
&=& ( \dfrac{n+2}{2n(n+1)} + \dfrac{1}{2(n+1)} ) \til{R} \til{h}_{\al \obe}\, .
\end{eqnarray*}

\n Thus, $$\til{R}_{\al \obe} = \frac{\til{R}}{n} \til{h}_{\al \obe}\, .$$
\end{proof}

\section{An ellipticity property}

In this section, we provide a short proof of Theorem 1.5. We will base on the following lemmas whose proof can be found in \cite{CNS}:

\begin{lema}
\label{Lema}
\n For any two hermitian matrices $A, B \in \Gamma_{k}^+$ of order $n$ and any $t
\in [0,1]$, we have the following inequality

\[
\{\si_k((1-t)A +tB)\}^{\frac{1}{k}} \geq (1-t)\{\si_k(A)\}^{\frac{1}{k}}+t\{\si_k(B)\}^{\frac{1}{k}}
\]

\n Furthermore, if $A \in \Gamma_{k}^+$, then $T_{k-1}(A)$ is positive definite.
\end{lema}

\begin{lema}
\label{lema4} Let $\La=(\la_1, \ldots , \la_n) \in \Ga^{+}_{k}$. Then,

\[
\si_{k-1}(\La) \geq \frac{k}{n-k+1} C(n,k)^{1/k}[\si_k(\La)]^{(k-1)/k}\, .
\]

\end{lema}

With these lemmas at hand, we prove Theorem 1.5 below.

\begin{proof}
If $\ta$ is $k$-positive, we have, by definition, that $S_{\ta} \in \Ga_{k}^{+}$. Since $M$ is compact, at a maximum point $p$ of the
function $u$, we obtain

\[
S_{\tta}(p) = S_{\ta}(p) -2\grad^{2}u(p)
\]

\n where $\grad^{2}u(p)$ is negative semi-definite. By Lemma \ref{Lema}, it follows that $S_{\tta}(p) \in \Ga_{k}^{+}$. So, using the assumption

\[
\si_k(S_{\tta}) = \mathcal L u > 0\ \ \mbox{on}\ \ M\, ,
\]

\n from the connectedness of the G\"{a}rding cones
$\Ga_{k}^{+}$ and from Lemma \ref{lema4}, we deduce that $S_{\tta} \in \Ga_{k}^{+}$.

\n The linearization of 

\[
F[u,\grad u, \grad^2 u] = \si_k({ S_{\ta} -2\grad^{2}u+(iT_{\ta}u-\normc {du} 2 \ta )L_{\ta} }) \, .
\]

\n toward $\phi$ yields

\[
D F[u,\grad u, \grad^2 u](\phi)=\si_1\lp{ \til T_{k-1}(S_{\tta}) \circ \til \grad^2 \phi }\rp.
\]

\n The conclusion then follows from the second part of Lemma \ref{Lema}, since $\til T_{k-1}(S_{\tta})$ is positive definite.
\end{proof}

\section{The CR $k$-Yamabe constant}

The canonical structure $\hat{\ta}$ on the CR sphere $\s^{2n+1}$ is pseudo-Einstein and torsion-free. In particular, if $\tta \in \mathcal{C}[\hat{\ta}]_{+}$ has constant pseudohermitian $k$-curvature, then $\tta$ is pseudo-Einstein. So, $\tta$ has constant Webster scalar curvature. In \cite{JL3}, Jerison and Lee characterized these pseudohermitian structures. Moreover, they established that if $\tta$ has constant Webster scalar curvature, then $\tta$ is pseudo-Einstein and torsion-free. In this case, $\tta$ has constant pseudohermitian $k$-curvature. Corollary 1.4 then follows from Corollary 1.3 and from the characterization provided in Theorem A of \cite{JL3}.

On the other hand, the following question remains open: is there any pseudohermitian structure on $\s^{2n+1}$ compatible with the canonical structure which has constant pseudohermitian $k$-curvature and is not CR equivalent to the standard one? We believe not and raise the following conjecture:

\begin{conj}
Let $(\s^{2n+1}, \hat{\ta})$ be the CR sphere of dimension $2n + 1$. Then, $\tta \in [\hat{\ta}]_{+}$ has constant pseudohermitian
$k$-curvature if, and only if, it has constant Webster scalar curvature. In particular, any solution of the $k$-Yamabe problem on $\s^{2n+1}$ can be built from a scalar multiple of the canonical structure $\hat{\ta}$ by a CR-automorphism of the CR sphere.
\end{conj}

Our main target in this section is proving Theorems 1.3 and 1.4. Theorem 1.4 furnishes a partial answer to Conjecture 6.1.

Let $(M, \ta)$ be a CR manifold of dimension $2n + 1$. Given a conformal pseudohermitian structure $\tta = e^{2u}\ta$, we recall that

\[
S_{\tta} = S_\ta - 2\grad^2 u +(i T_\ta u - \|du\|^{2}_{\ta})L_{\ta}\, .
\]

\n Making the standard change $e^{2u}=v^{2/n}$, we get

\[
\grad^{2}u (T_{\al}, T_{\obe})=u_{\al \obe} = \frac{1}{n}( \frac{v_{\al \obe}}{v} - \frac{v_{\al} v_{\obe}}{v^2} ) \, ,
\]

\[
T_{\ta} u =\frac{1}{n} \frac{T_{\ta} v}{v}\, ,
\]

\[
\|du\|^{2}_{\ta} = L^{*}_{\ta}(du , du) = L^{*}_{\ta} (\frac{dv}{nv} , \frac{dv}{nv})=\frac{1}{n^2 v^2} \|dv\|^{2}_{\ta}
\, .
\]

\n So,

\[
S_{\tta} = S_{\ta} -\frac{2}{n}\frac{\grad^2 v}{v} +\frac{2}{n}\frac{\grad v \otimes \grad v}{v^2} + \frac{2}{n} (
\frac{i}{2} \frac{T_{\ta} v}{v} - \frac{\|dv\|^{2}_{\ta}}{2nv^2} ) L_{\ta}\, .
\]

\n But, this equality can be rewritten as

\[
S_{\tta} = \dfrac{2}{n} \dfrac{V}{v}\, ,
\]

\n where

\begin{equation} \label{3}
V = V[v] := \frac{nv}{2}S_{\ta} -\grad^2 v+\frac{\grad v \otimes \grad v}{v}+\frac{1}{2} ( i T_{\ta} v- \frac{\|dv\|^{2}_{\ta}}{nv} )
L_{\ta}\, .
\end{equation}

\n Thus,

\[
\si_k(\tta) = \si_k ( v^{-2/n}h^{\ga \obe} \dfrac{2}{n} \dfrac{V_{\al \obe}}{v} ) = {( \frac{2}{n} )}^k
v^{-k(1+2/n)} \si_k(V_{\al}^{\;\; \ga}) = {( \frac{2}{n} )}^k v^{-k(1+2/n)} \si_k(V)\, ,
\]

\n and so, if $\si_k(\tta)$ is a constant $c$, then $v$ satisfies

\begin{equation}
\label{opYam1} v^{(1-k)\tfrac{n+2}{n}}\si_k(V)=\la v^{\tfrac{n+2}{n}}=\la v^{p-1}
\end{equation}

\n where $\la = \frac{n^k}{2^k}k c$.

As in the introduction, given a compact CR manifold $(M, \ta)$ of dimension $2n + 1$, we introduce the constants

\[
\la_k(M) = \inf \{  Y_k(\tta) : \tta \in {\cal M} \cap [\ta]\}\, ,
\]

\[
\la_k^{+}(M) = \inf \{  Y_k(\tta) : \tta \in {\cal M} \cap [\ta]_{+}\}
\]

\n and

\[
\la_k^{\mathcal{C}}(M) = \inf \{  Y_k(\tta) : \tta \in {\cal M} \cap \mathcal{C}[\ta]_{+}\}\, .
\]

\n As a direct consequence of Theorem 1.1 and Corollary 1.4, we
derive the following result:

\begin{coro}
\label{cor5.1.3} Let $(\s^{2n+1}, \hat{\ta})$ be the CR sphere of dimension $2n + 1$. Then, the constants $\la_{k}(\s^{2n+1})$, $\la_{k}^+(\s^{2n+1})$ and $\la_{k}^{C}(\s^{2n+1})$ are achieved by multiples of $\hat{\ta}$ by some conformal CR-automorphism of $\s^{2n+1}$. Moreover, $\la_{k}^{C}(\s^{2n+1})$ is achieved only by these structures.
\end{coro}

For $\til{\ta} = v^{p - 2} \ta$, we already know that

\[
\si_k(\tta) = {( \frac{2}{n} )}^k v^{k(1 - p)} \si_k(V)
\]

\n and

\[
dV_{\tta} = v^{p}\; dV_{\ta}\, ,
\]

\n where $p = 2 + 2/n$ and $V$ is provided in (\ref{3}).

We now turn our attention to the normalized functional

\[
J_k(v) := \frac{Y_k(\tta)}{V_{\tta}(M)^{1 - 2k/np}}= \frac{\int_M v^{p(1-k)+k}\si_k \lp {V}\rp\; dV_{\ta}}{{\lp \int_M v^p\; dV_{\ta}\rp}^{1 - 2k/np}}\, .
\]

\n This will be called the \textbf{CR $k$-Yamabe functional}. When $k = 1$, clearly $J_k$ is, modulo a constant factor, the CR Yamabe functional (see \cite{JL4}).

Notice that

\[
\la_k(M) = \inf \{ {( \frac{2}{n} )}^k J_k(v) : v \in C^{\infty}(M) \; , \; v > 0\ \ \mbox{on} \ M\}\, .
\]

\n So, if $\tta = v^{p-2} \ta$ satisfies $C_{\tta}=0$, we conclude
that $\tta$ has constant pseudohermitian $k$-curvature if, and only
if, $v$ is a critical point of the CR $k$-Yamabe functional.

We next furnish the proof of Theorem 1.3.

\begin{proof}

We first recall that

\[
\s^{2n+1}_{*}:=\s^{2n+1}-\{(0,-1)\}
\]

\n can be identified, via the Cayley transformation $F:S^{2n+1}_{*} \rightarrow \He^n$, with the Heisenberg group
$\He^n$. Moreover, we can consider the pseudohermitian structure $\ta_0$ on $\s^{2n+1}_{*}$ inherited from the canonical pseudohermitian structure $\Ta_0$ on $\He^n$, namely $\ta_0 = F^* \Ta_0$. Thus, the map

\[
F:(\s^{2n+1}_{*}, \ta_0) \rightarrow (\He^n , \Ta_0)
\]

\n is an isopseudohermitian map.

\n The CR $k$-Yamabe constant on $\s^{2n+1}$ can be written as

\[
\la_k(\s^{2n+1})=\inf \{ { ( \frac{2}{n} )^k \frac{\int_{\He^n} v^{p(1-k)+k}\si_k \lp {V}\rp\; dV_{\Ta_0}}{{\lp \int_{\He^n} v^p\; dV_{\Ta_0}\rp}^{1 - 2k/np}} : v \in C^{\infty}(\He^n) \; , \; v>0\ \ \mbox{on} \ M}\}\, .
\]

\n Consider now the function $v_0 \in C^{\infty}(\He^n)$ given by

\[
v_0 = v_0(z,t)=|w+i|^{-n}\, ,
\]

\n where $w=t+i|z|^2$. We can write

\[
v_0^{2/n} \Ta_0 = |w+i|^{-2}\Ta_0\, .
\]

\n Letting

\[
b(z,t)=|w+i|^{-2}\, ,
\]

\n we have

\[
F^*(b \Ta_0)=(b \circ F)  F^* \Ta_0\, .
\]

\n On the other hand,

\begin{eqnarray*}
b \circ F (z,z^{n+1}) &=& b (\frac{z}{1+z^{n+1}} , -i \frac{z^{n+1} - \overline{z}^{n+1}}{|1+z^{n+1}|^2})\\
&=&|1+z^{n+1}|^2
\end{eqnarray*}

\n and it is well-known that

\[
\ta_0 = F^* \Ta_0 = |1+z^{n+1}|^{-2} \hat{\ta}\, .
\]

\noindent Thus,

\begin{eqnarray*}
v_0^{2/n} \ta_0 &=& F^* (v_0^{2/n} \Ta_0) = F^* (|w+i|^{-2} \Ta_0)= \hat{\ta}\, .
\end{eqnarray*}

\n Since $\hat{\ta}$ has constant pseudohermitian $k$-curvature and
is torsion-free, the function $v_0 \in C^\infty(\He)$ is an extremal
for the CR $k$-Yamabe functional

\[
J_k(v) = \frac{\int_{\He^n} v^{p(1-k)+k}\si_k ( {V})\; dV_{\Ta_0}}{{( \int_{\He^n} v^p\; dV_{\Ta_0})}^{1 - 2k/np}}\, .
\]

\n In other words,

\[
\la_k(\s^{2n+1})={( \frac{2}{n} )}^k \frac{\int_{\He^n} v_0^{p(1-k)+k}\si_k ( {V_0})\; dV_{\Ta_0}}{{( \int_{\He^n} v_0^p\; dV_{\Ta_0})}^{1 - 2k/np}}\, ,
\]

\n where

\[
V_0 =\frac{nv_0}{2}S_{\Ta_0} -\grad^2 v_0 + \frac{\grad v_0 \otimes \grad v_0}{v_0} + \frac{1}{2} ( i T_{\Ta_0} v_0- \frac{\|dv_0\|^{2}_{\Ta_0}}{nv_0} )
L_{\Ta_0}\, .
\]

\n But,

\begin{eqnarray*}
\la_k(\s^{2n+1})&=& \frac{\int_{\He^n} ( \frac{2}{n} )^k v_0^{p(1-k)+k}\si_k ( {V_0})\; dV_{\Ta_0}}{{( \int_{\He^n} v_0^p\; dV_{\Ta_0})}^{1 - 2k/np}}\\
&=& \dfrac{\int_{\He^n} \si_k (v_0^{2/n}\Ta_0)\; dV_{v_0^{2/n}\Ta_0}}{{( \int_{\He^n}\; dV_{v_0^{2/n}\Ta_0})}^{1-2k/np}}\\
&=& \dfrac{\int_{\s^{2n+1}} \si_k (\hat{\ta})\; dV_{\hat{\ta}}}{{( \int_{\s^{2n+1}}\; dV_{\hat{\ta}})}^{1-2k/np}}\, .
\end{eqnarray*}

\n Since $\si_k (\hat{\ta})$ is constant, we derive

\[
\la_k(\s^{2n+1})= \si_k (\hat{\ta}) \frac{\int_{\s^{2n+1}}\; dV_{\hat{\ta}}}{{( \int_{\s^{2n+1}}\; dV_{\hat{\ta}})}^{1-2k/np}} =
\si_k (\hat{\ta}) {\lp \int_{\s^{2n+1}}\; dV_{\hat{\ta}}\rp}^{2k/np}\, .
\]

\n Using the fact that $\hat{\ta}$ is pseudo-Einstein, we have

\[
\si_k(\hat{\ta}) = C(n,k) \frac{{\hat{R}}^k}{{\lp 2n(n+1) \rp}^k} \, ,
\]

\n where $\hat{R}$ denotes the Webster scalar curvature of $\hat{\ta}$.

\n Therefore,

\[
\la_k (\s^{2n+1}) = \frac{C(n,k)}{{\lp 2n(n+1) \rp}^k} {\lp \hat{R}
{\lp  \int_{\s^{2n+1}}\; dV_{\hat{\ta}}\rp}^{2/np} \rp}^k\, .
\]

\n Thanks to the fact of the Webster scalar curvature $\hat{R}$ to be constant, the above equality yields

\[
\la_k (\s^{2n+1}) = \frac{C(n,k)}{{\lp 2n(n+1) \rp}^k} {\lp \frac{\int_{\s^{2n+1}} \hat{R}\; dV_{\hat{\ta}}}{{\lp  \int_{\s^{2n+1}}\;
dV_{\hat{\ta}}\rp}^{2/p}} \rp}^k\, .
\]

\n Finally, consider the CR Yamabe functional

\[
I(\ta) = \frac{\int_{\s^{2n+1}} R\; dV_{\ta}}{{( \int_{\s^{2n+1}}\; dV_{\ta})}^{2/p}}
\]

\n defined on $[\hat{\ta}]_+$, where $R$ denotes the Webster scalar curvature of $\ta$. In \cite{JL3}, Jerison and Lee proved that $\hat{\ta}$ is an extremal of this functional and that

\[
\la(\s^{2n+1})=\inf \{ I(\theta) : \theta \in {\cal M} \cap [\hat{\ta}]_+\} = \frac{\int_{\s^{2n+1}} \hat{R}\; dV_{\hat{\ta}}}{{(  \int_{\s^{2n+1}}\; dV_{\hat{\ta}})}^{2/p}} = 2\pi n(n+1)\, .
\]

\n Consequently,

\[
\la_k (\s^{2n+1})=\frac{C(n,k)}{{\lp 2n(n+1) \rp}^k} {{\lp 2\pi
n(n+1) \rp}^k}=C(n,k) \pi^k\, .
\]

\end{proof}

We conclude this section with the proof of Theorem 1.4.

\begin{proof}
Let $\tta=v^{p - 2} \ta$ be a $k$-positive pseudohermitian structure
conformal to $\ta$. Since $\si_1(\tta)$ is positive, the Webster
scalar curvature of $\tta$ is positive. We have the operator
$\si_1(\la(S_{\ta}))$ satisfies the comparison principle, see
\cite{bony} and \cite{LiMo}.

Let $v_1$ be a solution of the CR Yamabe problem. As it is well-known,

\[
J_1(v_1) \leq \la(\s^{2n+1})\, .
\]

\n Moreover, when $(M,\ta)$ is non-locally CR-equivalent to
$\s^{2n+1}$, we have the strict inequality $J_1(v_1) <
\la(\s^{2n+1})$. By Theorem 1.5, proved in the previous section, and
the standard elliptic theory developed in \cite{CNS}, \cite{EL2},
\cite{EL3}, \cite{EL4} and \cite{EL5} (see also \cite{EL6}), the
equation

\[
\si_k(V)=C_{n,\, k} v_1^{(p-1)k}\ \ \mbox{on}\ \ M\, ,
\]

\n admits a unique positive solution $v_k \in C^\infty(M)$, where $C_{n,k} = C(n,k)n^k$. By Lemma \ref{lema4}, we obtain

\[
- p \lap_b v_k + Rv_k=\si_1\lp \la(V_k)\rp \geq n^2 v_{1}^{p-1}\ \ \mbox{on}\ \ M\, .
\]

\n On the other hand, $v_1$ satisfies

\[
-p \lap_b v_1 + R v_1 = n^2 v_{1}^{p-1}\ \ \mbox{on}\ \ M\, .
\]

\n So, by the comparison principle,

\[
v_k \geq v_1\ \ \mbox{on}\ \ M\, .
\]

\n In order to establish that

\[
\la_k^{+}(M) \leq \la_k^{+}(\s^{2n+1})\, ,
\]

\n it suffices to show that

\[
J_k(v_k) \leq \la^{+}_k(\s^{2n+1})\, .
\]

\n Since, for any $k \geq 2$,

\[
p(1-k) + k = p-k(p-1) < 0\, ,
\]

\n we have

\[
v_k^{p(1-k)+k} \leq v_1^{p(1-k)+k}\ \ \mbox{on}\ \ M\, .
\]

\n Thus,

\begin{eqnarray*}
\int_M v_k^{p(1-k)+k}\si_k \lp {V_k}\rp\; dV_{\ta} &\leq& C_{n,k} \int_M v_1^{p(1-k)+k} v_1^{(p-1)k}\; dV_{\ta}\\
&\leq& C_{n,k} \int_M v_1^{p(1-k)+k}v_1^{(p-1)k}\; dV_{\ta}\\
&=& C_{n,k} \int_M v_1^{p}\; dV_{\ta}\, .
\end{eqnarray*}

\n But, we also have

\[
\int_M v_k^{p}\; dV_{\ta} \geq \int_M v_1^{p}\; dV_{\ta}\, .
\]

\n Therefore,

\[
J_k(v_k) = \frac{\int_M v_k^{p(1-k)+k} \si_k \lp {V_k}\rp\; dV_{\ta}}{{\lp \int_M v_k^p\; dV_{\ta}\rp}^{1 - 2k/np}} \leq C_{n,k} \lp \int_M v_1^{p}\; dV_{\ta}\rp^{2k/np}
\]

\n and so, the inequality $J_k(v_k) \leq \la^{+}_k(\s^{2n+1})$ follows directly from

\[
\la^{+}_k(\s^{2n+1}) = \frac{C(n,k)}{2n(n+1)}{ \la(\s^{2n+1}) }^k\, .
\]

\end{proof}

\n {\bf Acknowledgments:} The first and third authors was partially
supported by the CNPq and Fapemig financial support agencies.

\end{document}